\documentclass[aip,jcp,showpacs,showkeys,a4paper,10pt]{revtex4-1}

\usepackage{graphicx}
\usepackage{amsmath}
\usepackage{amsfonts}
\usepackage{amssymb}
\usepackage{multirow}
\usepackage{rotating}
\usepackage{float}
\usepackage{hyperref}

\begin{document}
 
 \title{Influence of mass and potential energy surface geometry on roaming in Chesnavich's CH$_4^+$ model} 
 \begin{abstract}
  Chesnavich's model Hamiltonian for the reaction CH$_4^+ \rightarrow$ CH$_3^+$ + H is known to exhibit a range of interesting dynamical phenomena including roaming. The model system consists of two parts: a rigid, symmetric top representing the $\text{CH}_3^{+}$ ion, and a free H atom. We study roaming in this model with focus on the evolution of geometrical features of the invariant manifolds in phase space that govern roaming under variations of the mass of the free atom $m$ and a parameter $a$ that couples radial and angular motion. In addition we establish an upper bound on the prominence of roaming in Chesnavich's model. The bound highlights the intricacy of roaming as a type of dynamics on the verge between isomerisation and nonreactivity as it relies on a generous access to the potential wells to allow reactions as well as a prominent area of high potential that aids sufficient transfer of energy between the degrees of freedom to prevent isomerisation.
 \end{abstract}
 \author{Vladim{\'i}r Kraj{\v{n}}{\'a}k}
 \affiliation{School of Mathematics, University of Bristol, Bristol, BS8 1TW, UK}
 \author{Stephen Wiggins}
 \affiliation{School of Mathematics, University of Bristol, Bristol, BS8 1TW, UK}
 \maketitle

 \section{Introduction}
  Following recent developments in the understanding of the roaming mechanism in molecular dynamics \cite{Townsend2004,Bowman2011,Bowman2014,Mauguiere2014,Mauguiere2014b,Mauguiere2015,Mauguiere2016,Bowman2017,Krajnak2018}, the results of Heazlewood et al. \cite{Heazlewood08} on roaming being the dominant mechanism for acetaldehyde photodissociation triggered the following question: Can some of the standard parameter values in Chesnavich's CH$_4^+$ model \cite{Chesnavich1986} be altered so that it admits roaming as the dominant form of dissociation? It is well known that roaming plays a minor role in the dissociation of Formaldehyde\cite{Townsend2004}. Since acetaldehyde dissociates into CH$_4+$CO, one of the obvious differences in the dissociation of Formaldehyde into H$_2+$CO is the mass ratio of the products. Therefore we investigate the influence of various masses of the free atom on roaming in Chesnavich's CH$_4^+$ model. Furthermore we consider different values of parameter $a$, the coupling parameter of the angular and radial degrees of freedom in Chesnavich's potential, and deduce its impact on roaming. Physically speaking, $a$ controls the strength of the short range anisotropy of the CH$_3^+$ molecule and its precise role in the potential is discussed in Section \ref{sec:potential}.
 
  Chesnavich's CH$_4^+$ model \cite{Chesnavich1986} is a phenomenological $2$ degree of freedom model introduced that was introduced to study the transition from vibrational reactants to rotational products in the presence of multiple transition states. We use it to study the first stage of roaming, where a hydrogen atom separates from the rigid CH$_3^+$ molecule and instead of dissociating, it roams in a region of nearly constant potential only to return to the molecule without forming a bond. The process of intramolecular abstraction and subsequent dissociation requires both H atoms to be mobile and thus requires at least $4$ degrees of freedom. At the moment no tools for a qualitative dynamical study that would explain roaming in its entirety in $4$ degrees of freedom have been developed.
  
%  \section{Chesnavich's CH$_4^+$ model}
  We study the system in a polar centre of mass frame, where $r$ is the distance of the mobile H atom from the centre of mass of the system in \AA$ $ and $\theta$ is the angle describing the relative orientation of H and CH$_3^+$ in radians. The momenta $p_r$ and $p_\theta$ are cannonically conjugate to $r$ and $\theta$ respectively.
  
  The model we study is defined by the Hamiltonian
  $$H(r,\theta,p_r, p_\theta) = \frac{1}{2 \mu} p_r^2 + \frac{1}{2}p_\theta^2 \left(\frac{1}{I} + \frac{1}{\mu r^2}\right)  + U(r,\theta),$$
  where $\mu=\frac{m_{CH_3}m_{H}}{m_{CH_3}+m_{H}}$, with $m_{H}=1.007825$ u and $m_{CH_3}=12.0$ u, is the reduced mass of the system, and $I=2.373409$~u\AA$^2$ is the moment of inertia of the rigid body CH$_3^+$. Chesnavich's potential $U(r,\theta)$ has the form
  \begin{equation}\label{eq:U}
   U(r,\theta ) = U_{CH} (r) + U_{coup} (r,\theta),
  \end{equation}
  and it is the sum of radial long range potential $U_{CH}$ and short range ``hindered rotor'' potential $U_{coup}$, that represents the anisotropy \cite{Chesnavich1986,Jordan1991} of the rigid molecule CH$_3^+$. We use the standard definition of $U_{CH}$ and $U_{coup}$ and standard parameter values as suggested by Chesnavich \cite{Chesnavich1986} and used in recent publications \cite{Mauguiere2014,Mauguiere2014b,Krajnak2018}, namely
  \begin{equation}\label{eq:UCH}
   U_{CH} (r) =  \frac{D_e}{c_1 - 6} \left( 2 (3-c_2) e^{c_1 (1-x)}  - \left( 4 c_2 - c_1 c_2 + c_1 \right) x^{-6} - (c_1 - 6) c_2 x^{-4} \right), 
  \end{equation}
  where $x = \frac{r}{r_e}$, $D_e=47$ kcal/mol, $r_e=1.1$ \AA, $c_1=7.37$ and $c_2=1.61$, and
  \begin{equation}\label{eq:Ucoup}
   U_{coup} (r,\theta) = \frac{U_e e^{-a(r-r_e)^2}}{2} (1 - \cos 2 \theta ),
  \end{equation}
  where $U_e=55$ kcal/mol. The parameter $a$ (in \AA$^{-2}$) influences the value of $r$ at which the transition from vibration to rotation occurs, for example $a=1$ represents a late transition \cite{Mauguiere2014b,Mauguiere2014,Krajnak2018} and $a=4$ an early transition \cite{Mauguiere2014b}. In this paper we shall explore all values of $a$ that may be relevant to roaming.
  
  The total energy $H(r,\theta,p_r, p_\theta)=E$ is given in kcal/mol with respect to the dissociation energy $0$. 
  
  We introduce the potential in more detail and explain the role of the parameter $a$ is Section \ref{sec:potential}. Subsequently we give an overview of roaming up to date in Section \ref{sec:roam}, we explain how phase space structures influence the dynamics in Section \ref{sec:phase transport} and focus on how these structures cause roaming in \ref{sec:phase roaming}.
  
  To study the dependence of roaming on the mass of the free atom, we replace the free H atom by a atom of mass $m$, so that the reduced mass of the system is $\mu_m=\frac{m_{CH_3}m}{m_{CH_3}+m}$ and the Hamiltonian then is
  \begin{equation}\label{eq:Hamm}
   H_m(r,\theta,p_r, p_\theta) = \frac{1}{2 \mu_m} p_r^2 + \frac{1}{2}p_\theta^2 \left(\frac{1}{I} + \frac{1}{\mu_m r^2}\right)  + U(r,\theta).
  \end{equation}
  In this work we consider $m<m_{CH_3}$, since the free atom is usually the lighter of the two dissociated products. For comparison, the free atom would have to weight $m=8.57$ to be in the same proportion as in acetaldehyde.
  
  Since $m$ is only present in the kinetic part of the Hamiltonian (\ref{eq:Hamm}), its influence cannot be seen in configuration space. We explain how variations in mass of the free atom influence the system, phase space structures and roaming in different parts of phase space in Sections \ref{sec:inner}, \ref{sec:outer} and \ref{sec:small}.

  The phase space structures that enable us to study dynamics in phase space are normally hyperbolic invariant manifolds (NHIMs) and the corresponding stable and unstable invariant manifolds. We build upon the understanding of the role of NHIMs and invariant manifolds in governing dynamics as described in Refs. \onlinecite{OzoriodeAlmeida90, Waalkens10, Wiggins16} in great detail.
   
 \section{Role of $a$ in Chesnavich's CH$_4^+$ potential}\label{sec:potential}
  The potential $U(r,\theta)$ as defined by (\ref{eq:U}), (\ref{eq:UCH}) and (\ref{eq:Ucoup}), has the following characteristics:
  \begin{itemize}
   \item potential wells near $\theta=0$ and $\theta=\pi$ representing the two isomers of CH$_4^+$,
   \item areas of high potential near $\theta=\frac{\pi}{2}$ and $\theta=\frac{3\pi}{2}$ creating a potential barrier between the two wells,
   \item an area where the potential is monotonic and nearly constant due to $U(r,\theta)\in o(r^{-4})$ as $r\rightarrow\infty$, representing the dissociated state.
  \end{itemize}
  
  Note that due to the rotational symmetry $U(r,\theta)=U(r,\theta+\pi)$ and reflectional symmetry $U(r,\theta)=U(r,-\theta)$ of $U$ that follows from the anisotropic term $U_{coup}$ in (\ref{eq:Ucoup}), the wells and areas of high potential are symmetric.
  
  All of the characteristics listed above can be derived from the critical points of $U(r,\theta)$. Table \ref{table:equil} shows critical points of $U$ for $0\leq\theta<\pi$, while symmetric counterparts exist in $-\pi\leq\theta<0$. We denote critical points by $q$, the subscript indicates the index of the saddle and the superscript indicates the half plane in which the critical point lies: $+$ for $0\leq\theta<\pi$, $-$ for $-\pi\leq\theta<0$.
  
  \begin{table}[H]
    \begin{center}
    \begin{tabular}{c|c|c|c|c}
%     \hline\hline
    Energy (kcal mol$^{-1}$) & $r$ (\AA) & $\theta$ (radians) & Significance & Label \\
    \hline%\hline
    $-47$ & $1.1$ & $0$ & potential well & $q_0^+$ \\
    $>0$ & $<1.1$ & $\pi/2$ & isomerisation saddle & $q_1^+$ \\
    $>0$ & $>1.1$ & $\pi/2$ & local maximum & $q_2^+$ \\
    $<0$ & $>1.1$ & $\pi/2$ & isomerisation saddle & $\widetilde{q}_1^+$ \\
%     \hline\hline
    \end{tabular}
    \end{center}
    \caption{\label{table:equil} Equilibrium points for potential $U(r, \theta)$. Energy and radial coordinate of $q_1^+$, $q_2^+$ and $\widetilde{q}_1^+$ varies with $a$ and is shown graphically in Figure \ref{fig:thpi2}.}
   \end{table}
  
\begin{sidewaysfigure}
   \centering
   \includegraphics[width=.43\textwidth]{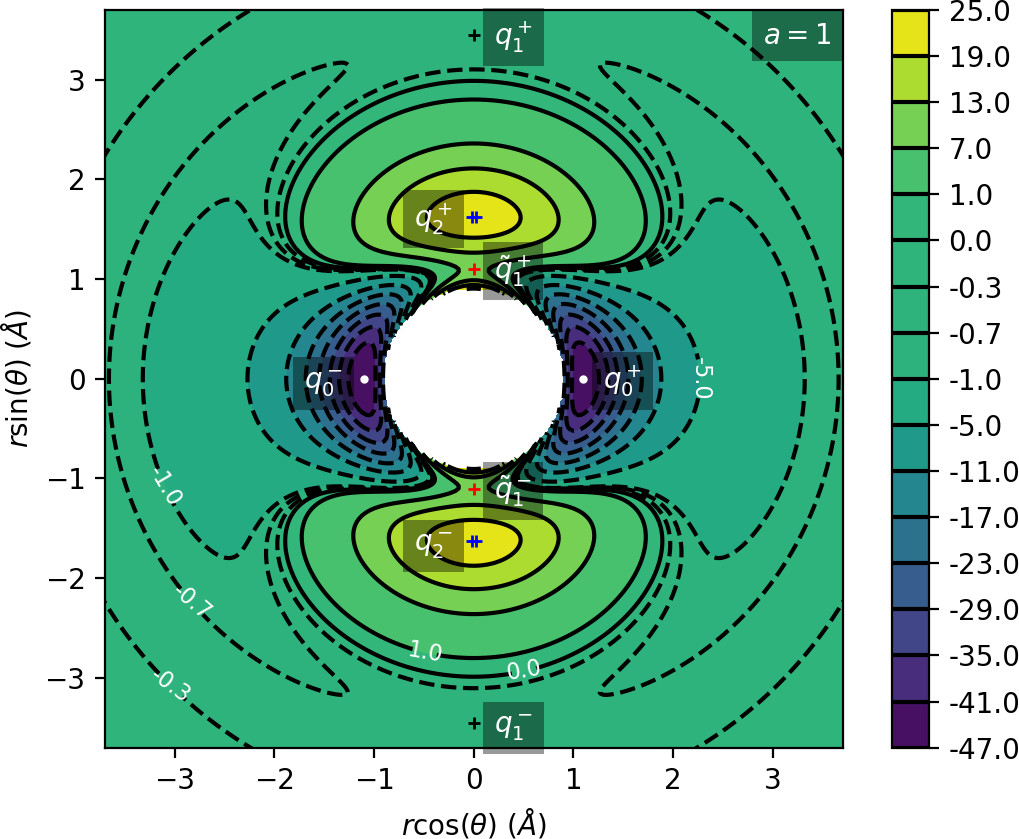}
   \includegraphics[width=.43\textwidth]{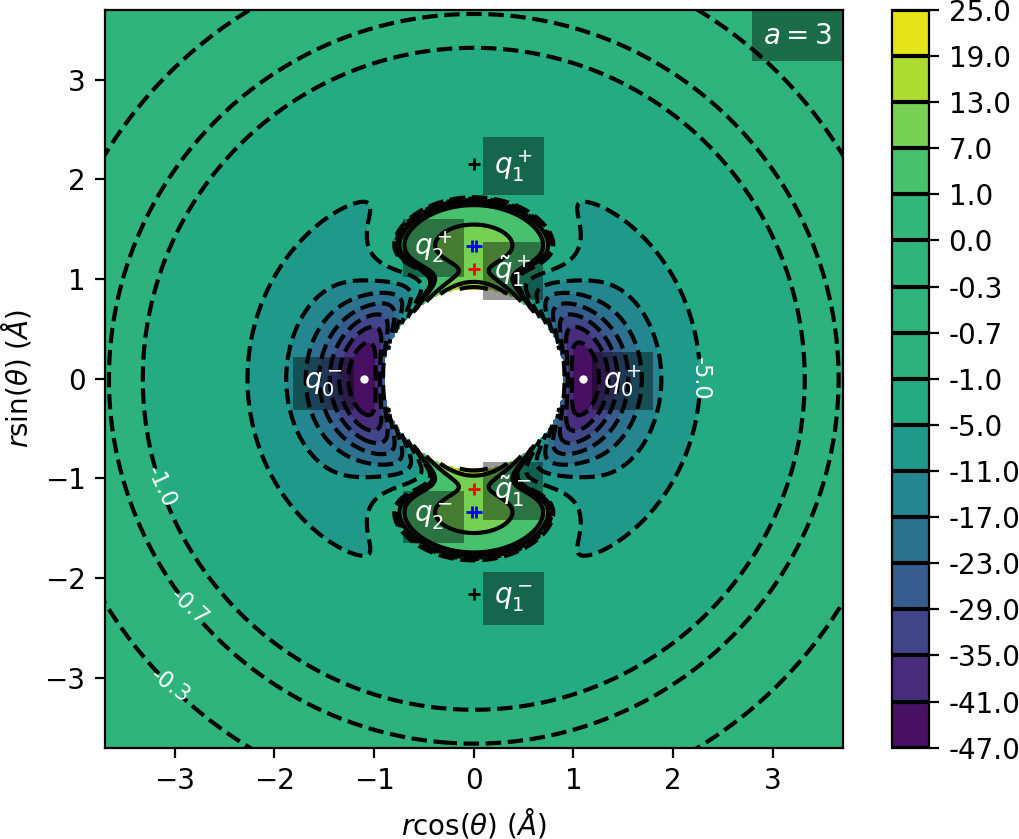}\\
   \includegraphics[width=.43\textwidth]{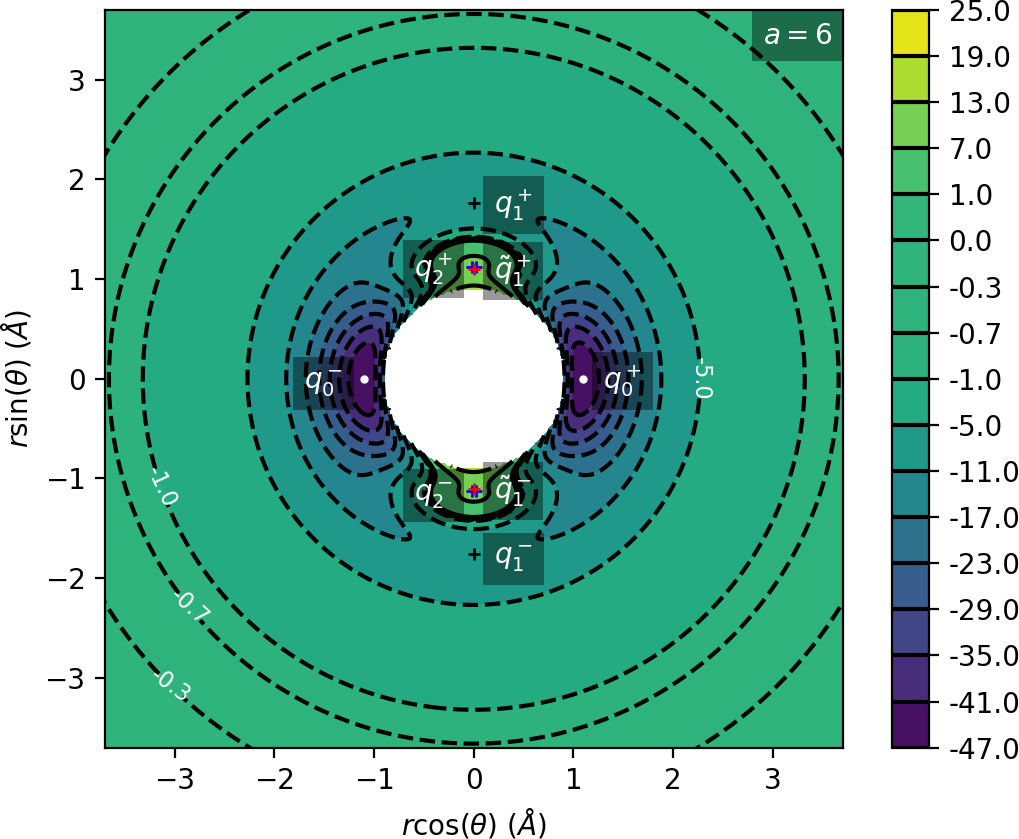}
   \includegraphics[width=.43\textwidth]{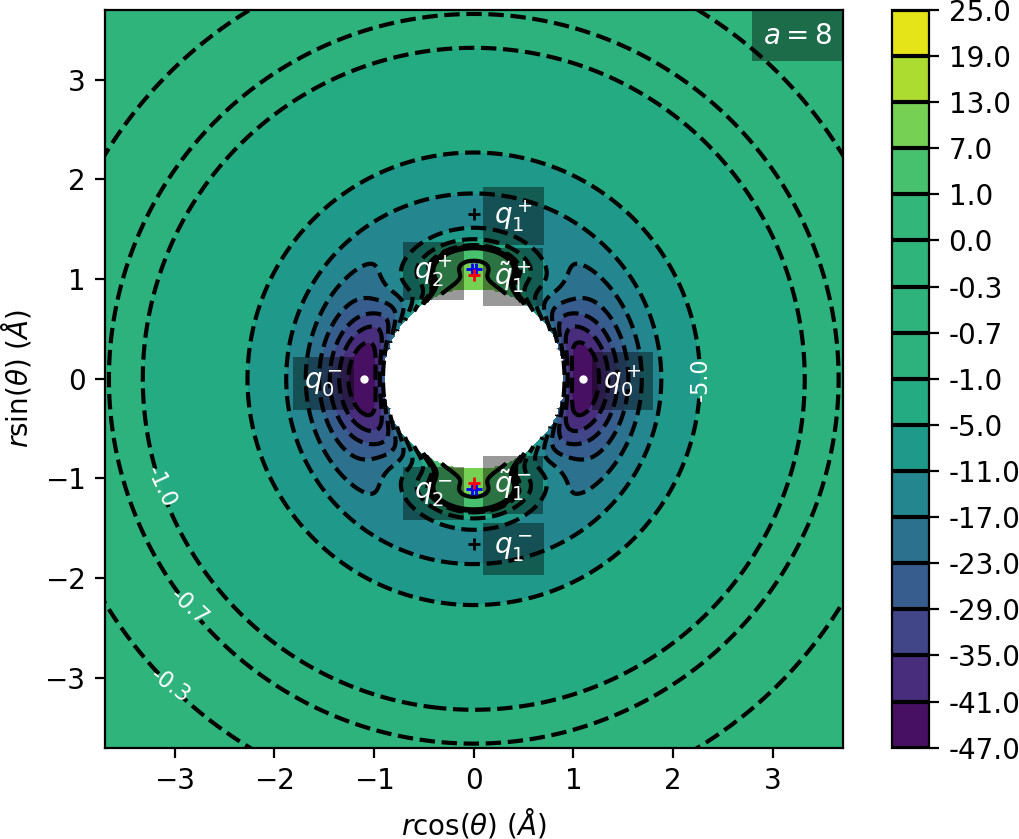}
   \caption{Contour plots of the potential energy surface $U$ for $a=1,3,6,8$. Dashed lines correspond to $U<0$, solid lines correspond to $U\geq0$. Contours correspond to values of potential shown on the colorbar right, with some values indicated in the plot. Critical points of the potential $q_0^-$, $q_1^-$, $q_2^-$ and $\widetilde{q}_1^-$, as introduced in Table \ref{table:equil} are indicated.}
   \label{fig:pot}
  \end{sidewaysfigure}

  We remark that $q_1^+$ and $q_2^+$ are energetically inaccessible at energies considered in this work. Four more critical points $q_0^-$, $q_1^-$, $q_2^-$ and $\widetilde{q}_1^-$ are related to the ones above by symmetry. The location of the critical points in configuration space can be seen in Figure \ref{fig:pot} for $a=1,3,6,8$.

\begin{figure}
   \centering
   \includegraphics[width=.9\textwidth]{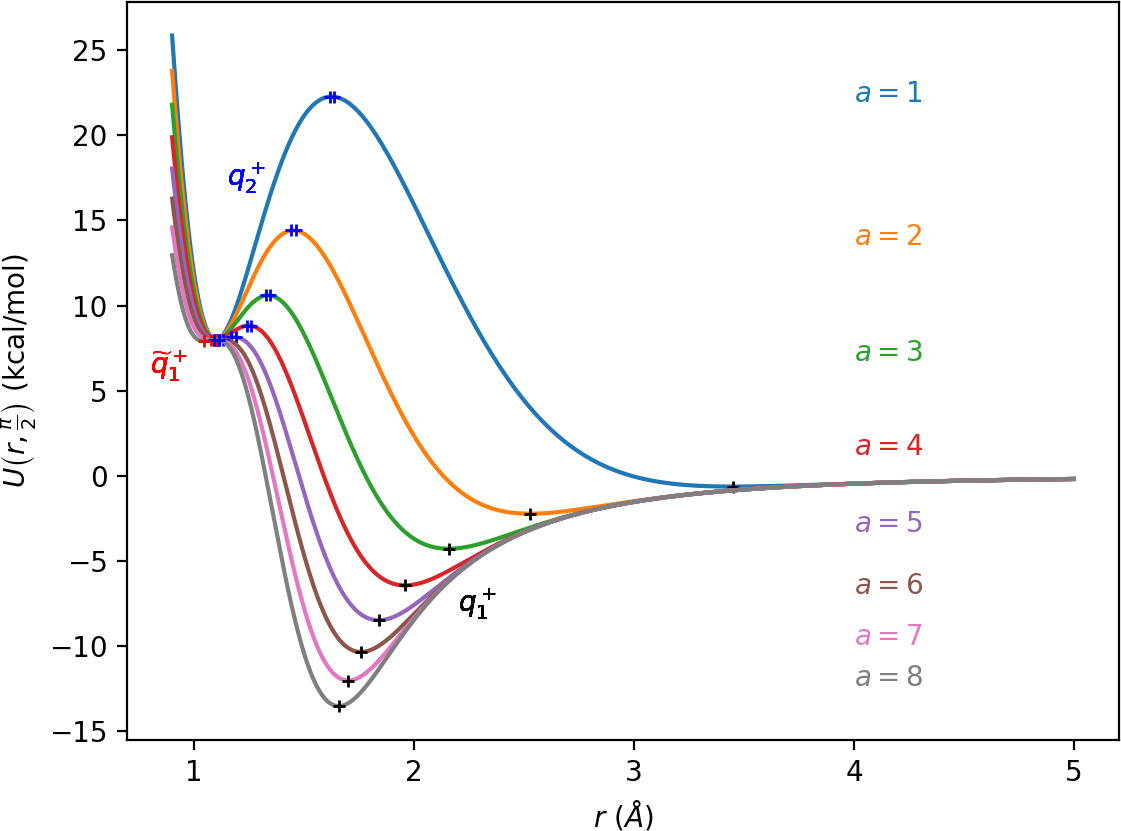}
   \caption{Radial sections of the potential $U$ along $\theta=\frac{\pi}{2}$ for $a=1,2,3,4,5,6,7,8$.}
   \label{fig:thpi2}
  \end{figure}

  The coupling term $U_{coup}$ through which $a$ influences the potential, vanishes around $\theta=0$ and $\theta=\pi$ and is maximal around $\theta=\pm\frac{\pi}{2}$. Therefore variations of $a$ leave most of the well unaffected, while the potential maxima and associated critical points vary significantly. Figure \ref{fig:thpi2} illustrates how the potential barrier between the wells recedes with increasing $a$. The figure also shows how $q_1^+$, $q_2^+$ and $\widetilde{q}_1^+$ are affected by $a$, note the decrease of energy of $q_1^+$. When considered in context of the contour plots in Figure \ref{fig:pot}, for $a=6$ and larger the wells merge into one from an energetic perspective, but dynamically remain distinct.

  Due to the exponential decay of the anisotropic term $U_{coup}$ and the $o(r^{-4})$ decay of $U$ as a whole, the area of flat potential that can be seen in Figures \ref{fig:pot} and \ref{fig:thpi2} remains unchanged. In this area $U$ is negative, monotonic in $r$, very close to zero and certainly has no saddles.
  
  Apart from different values of the potential, the wells also differ from the area of flat potential dynamically. In the potential well trajectories vibrate, i.e. oscillate in the angular direction, as they are bounded by the two areas of high potential. Outside of the wells, especially in the flat area, the absence of potential barriers enables unlimited movement in the angular direction that is manifested in trajectories completing full rotations around the centre of mass. The transition from one type of dynamics to the other depends heavily on the size and height of the areas of high potential that is influenced by $a$ as described above.
 
 \section{Roaming}\label{sec:roam}
  Roaming was discovered in the study of photodissociation of formaldehyde \cite{Townsend2004} (H$_2$CO) and explained the unusual energy distribution between H$_2$ and CO experimentally observed previously \cite{vanZee1993}. Since then roaming was reported in a number of other molecules \cite{Bowman2017}. Initial states of roaming resemble radical dissociation, when a single H atom escapes from HCO by breaking a single covalent bond and immediately dissociates. Instead of dissociating, the free hydrogen `roams' around a flat monotonic region of the potential near HCO and returns to the molecule to abstract the remaining hydrogen. Ultimately, the products of this dissociation are identical to the products of molecular dissociation. While molecular dissociation requires the system to pass over a potential saddle, no saddle is involved in roaming. Consequently the distribution of energy between the products can differ significantly. Following a study \cite{Bowman2011} of H$_2$CO and CH$_3$CHO it was established that regardless of similarity of products, roaming is actually closer to radical dissociation than to molecular dissociation and further evidence was published subsequently \cite{Huston2016}.
  
  In light of developments following a phase space approach to chemical reaction dynamics \cite{Wigginsetal01,Uzeretal02,Waalkens04}, we employ the definition of roaming introduced in Ref. \onlinecite{Mauguiere2014b}. This definition considers the number of intersections of a trajectory with a dynamically justified dividing surface and thereby accounts for the influence of various phase space structures. The dividing surface is constructed using an unstable periodic orbit that is not associated with any potential saddle. It was shown \cite{Mauguiere2015} how invariant manifolds of unstable periodic orbits convey trajectories between two potential wells in what is called a 'shepherding mechanism'. The exact phase space structure, an intersection of invariant manifolds of unstable periodic orbits, responsible for roaming was since identified \cite{Krajnak2018}. The key to understanding roaming follows from the use of toric surfaces of section \cite{MacKay2014, MacKay2015, Mauguiere2016} to study invariant manifolds.

 \section{Phase space structures governing transport}\label{sec:phase transport}
  There are three important families of periodic orbits \cite{Mauguiere2014,Mauguiere2014b,Krajnak2018} in Chesnavich's CH$_4^+$ model. By family of periodic orbits we mean a continuum of periodic orbits parametrised by energy, so that each family contains two orbits related by (reflectional) symmetry for a given $E$. We will refer to these orbits as the inner ($\Gamma^i$), middle ($\Gamma^a$) and outer ($\Gamma^o$) periodic orbit based on the radii of the orbits and their significance in dissociation. For a given $E$, the middle and outer families consist of two orbits with the same configuration space projections but opposite orientation, one with $p_\theta>0$ and one with $p_\theta<0$, hence clearly distinct in phase space. In a configuration space projection these orbits rotate counterclockwise and clockwise respectively. The inner family consists of two orbits, one on the edge of each potential well. In the following, statements regarding one of the orbits of a family automatically hold for the other one due to symmetry. If we only refer to one orbit, we refer to $\Gamma^a$ and $\Gamma^o$ with $p_\theta>0$ and $\Gamma^i$ oscillating around $\theta=0$.
  
\begin{figure}
   \centering
   \includegraphics[width=.49\textwidth]{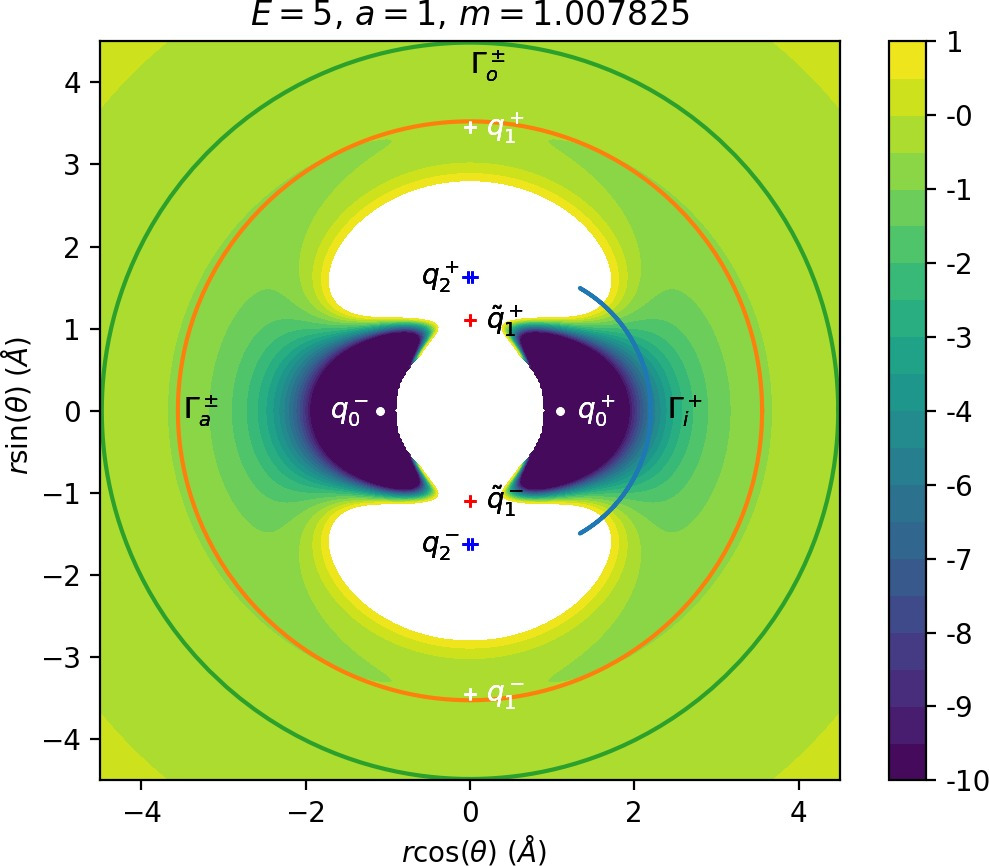}
   \includegraphics[width=.49\textwidth]{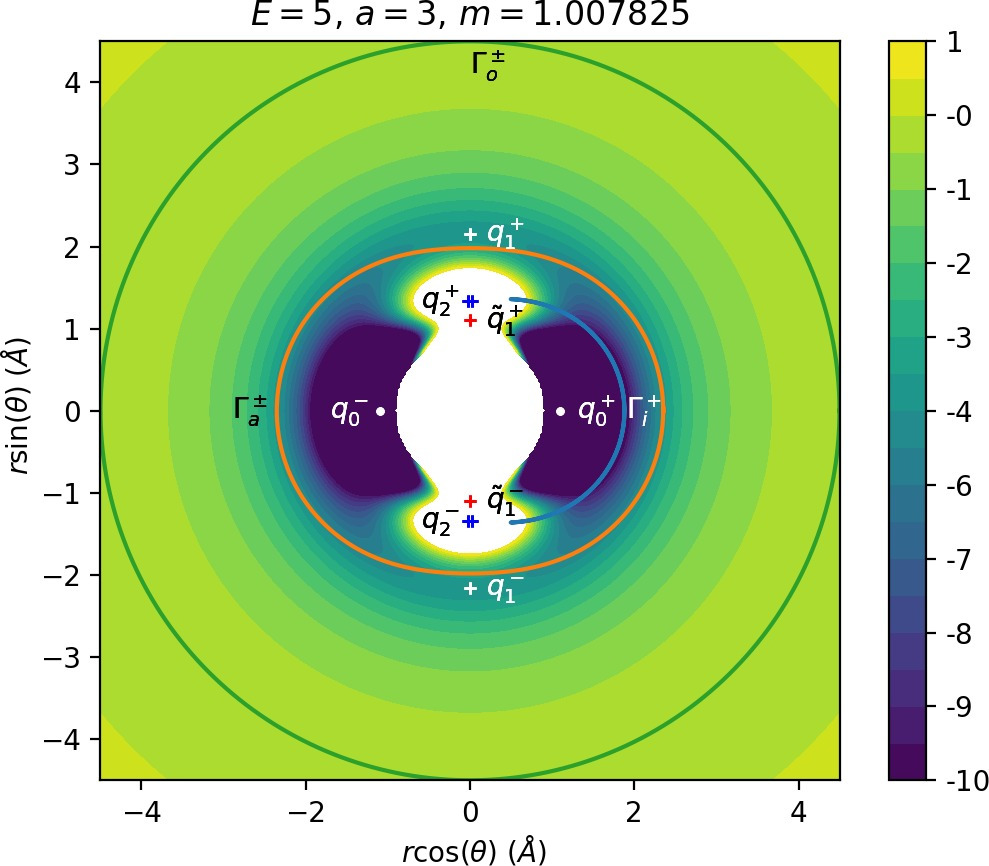}\\
   \includegraphics[width=.49\textwidth]{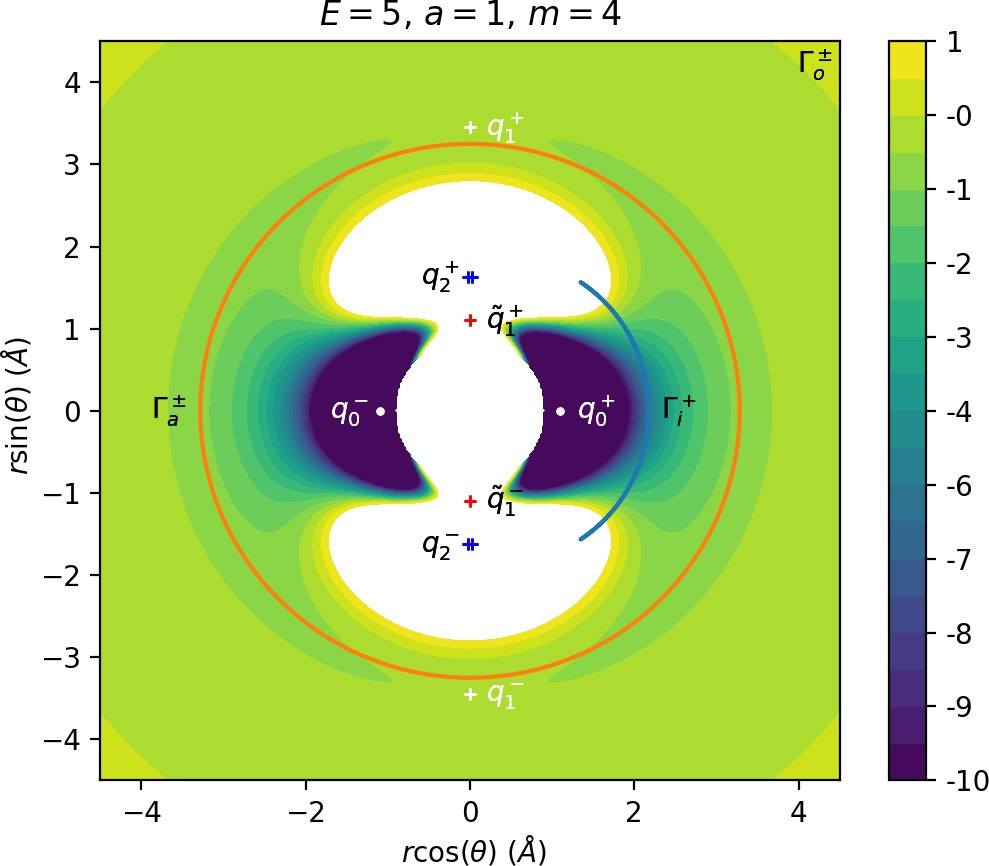}
   \includegraphics[width=.49\textwidth]{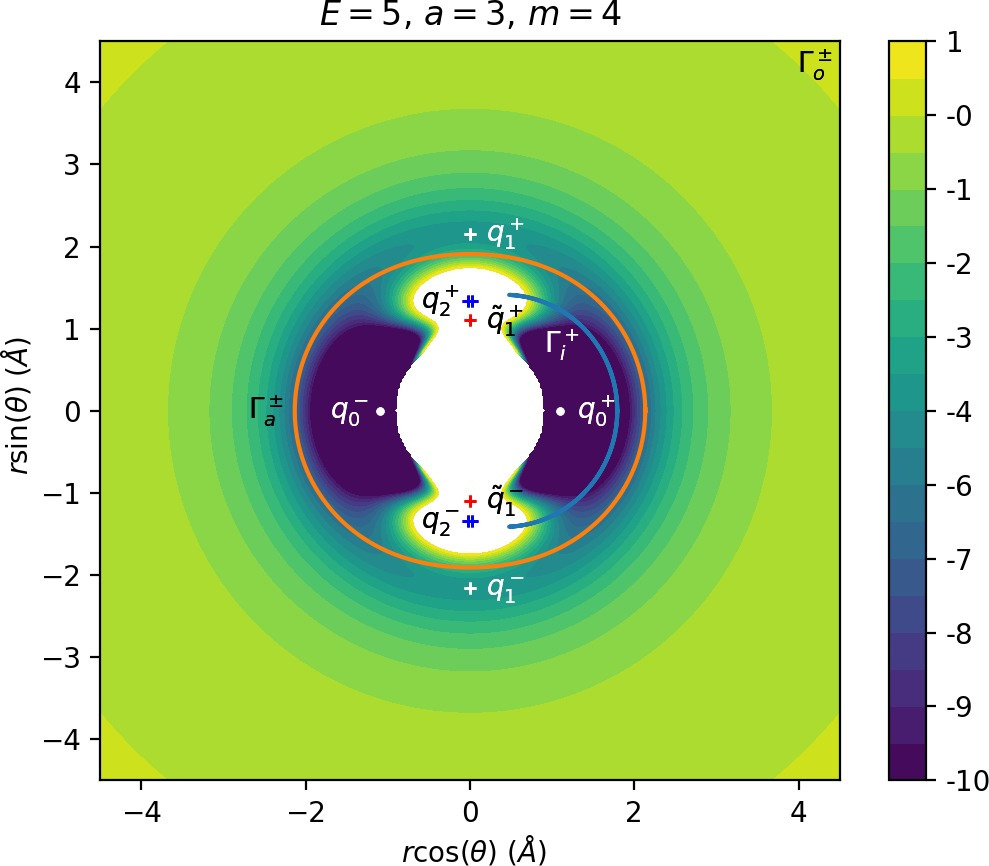}
   \caption{Configuration space projections of the inner, middle and outer periodic orbits for $E=5$ and a combination of $a=1,3$ and $m=m_H,4$. Note that for $m=4$, $\Gamma^o$ lies far outside the field of view. Critical points of the potential as introduced previously are also indicated.}
   \label{fig:poconfig}
  \end{figure}

  Configuration space projections are shown in Figure \ref{fig:poconfig}. Note that none of these families is associated with a saddle point on the potential energy surface.
  
  The family of outer periodic orbits $\Gamma^o$ is proven \cite{Krajnak2018} to exist due to a centrifugal barrier for a class of systems including the one considered here. It was also shown that the orbits in this family are unstable. States of the system beyond the outer orbits correspond to dissociated states CH$_3^++$H. Orbits of this family are therefore the outer boundary of the roaming region \cite{Mauguiere2014}. We prefer to use the more general term interaction region as the phase space structures responsible for roaming and isomerisation among trajectories that leave the well, also arbitrate whether incoming trajectories originating at $r=\infty$ react (reach the well) or not. In the context of transition state theory, this region is sometimes also called the collision complex.
  
  The potential wells represent the region where the CH$_4^+$ molecule is in a stable configuration and orbits of the inner family $\Gamma^i$ naturally delimit this region. All inner periodic orbits are unstable and form the inner boundary of the interaction region.
  
  Inside the interaction region lies the family of middle orbits, which are crucial in the definition of roaming. These orbits are unstable for small positive energies ($E<2.72$ for $a=1$ and $m=m_H$) and stable after a period doubling bifurcation involving a family denoted sometimes by FR$12$ \cite{Mauguiere2014,Mauguiere2014b} or $\Gamma^b$ \cite{Krajnak2018} that is of no further importance. The exact energy of the period doubling bifurcation varies with $m$ and $a$, but for $m=m_H$ and $a=1$ is at approximately $E=2.72$. The family is also subject to bifurcations at negative energies, but these are not relevant for roaming.
  
  The definition of roaming as introduced in Ref. \onlinecite{Mauguiere2014b} requires us to construct a dividing surface in phase space based on the configuration space projection of $\Gamma^a$. We shall denote the dividing surface associated with $\Gamma^a$ as DS$^a$. For a given energy $E$, DS$^a$ consists of all points $(r,\theta,p_r, p_\theta)$ on the energy surface $H_m(r,\theta,p_r, p_\theta)=E$, such that $(r,\theta)$ is a point on the configuration space projection of $\Gamma^a$.
  
  It is known that a spherical dividing surface may in the neighbourhood of an index-1 saddle bifurcate into a torus \cite{MacKay2014,MacKay2015} and it was recognised \cite{Mauguiere2016} that dividing surfaces constructed using periodic orbits that rotate such as DS$^a$, as opposed to vibrate, are tori. In addition, if $\Gamma^a$ is unstable, then DS$^a$ does not admit local recrossings.
  
  Indeed one can see that for every point $(r_{\Gamma^a},\theta_{\Gamma^a})$ on the configuration space projection of $\Gamma^a$, the corresponding phase space structure is a circle implicitly defined by (\ref{eq:Hamm}) as
  \begin{equation*}
   E-U(r_{\Gamma^a},\theta_{\Gamma^a}) = \frac{1}{2 \mu_m} p_r^2 + \frac{1}{2}p_\theta^2 \left(\frac{1}{I} + \frac{1}{\mu_m r^2}\right).
  \end{equation*}
  Note that the radius of the circle never vanishes along $\Gamma^a$ as the kinetic energy never vanishes along $\Gamma^a$.
  
  A trajectory is considered a roaming trajectory if it originates in the potential well and crosses DS$^a$ at least three times before dissociating. If a trajectory crosses DS$^a$ an even number of times and enters either one of the potential wells, it is an isomerisation trajectory.
  
  Roaming can be reformulated as a transport problem \cite{Krajnak2018} in the following manner. Let DS$^i$ and DS$^o$ be defined using $\Gamma^i$ and $\Gamma^o$ analogously to DS$^a$. It is well known that a dividing surface associated with a periodic orbit that reaches two equipotentials, such as $\Gamma^i$, is a sphere. Using reasoning similar to the above one can see, that because the kinetic energy along $\Gamma^i$ vanishes when it reaches the equipotentials, DS$^i$ consists of circles except for the poles which are points, thus a sphere. $\Gamma^i$, being the equator of this sphere, divides it into two hemispheres with unidirectional flux - trajectories leaving the well cross the outward hemisphere while trajectories entering the well cross the inward hemisphere.
  
  DS$^a$ and DS$^o$ are tori, but they too can be divided into two parts of unidirectional flux. This can be easily seen on DS$^o$, but a similar argument applies to DS$^a$. Since $\Gamma^o$ has a constant radius and $p_r=0$, the two orbits are characterised by $p_\theta^2$ being maximal for the given energy along the orbit. Any trajectory crossing DS$^o$ must therefore have $p_r\neq0$. All trajectories with $p_r>0$ cross the outward annulus of DS$^o$ bounded by the two orbits $\Gamma^o$, while all with $p_r<0$ cross the inward annulus and enter the interaction region.
  
  Note that because dividing surfaces divide the energy surface into two separate partitions, the dividing surfaces defined above have the following roles. Since $\Gamma^i$ lies on the edge of the potential well, the role of DS$^i$ is to divide the the energy surface region associated with the potential well from the rest of the energy surface. Every trajectory that enters or exists the well must cross DS$^i$. Similarly, DS$^o$ is used to divide the dissociated states from the rest of the energy surface and every trajectory going from one part of the energy surface into the other must cross DS$^o$. DS$^a$ divides the the energy surface interaction region between DS$^i$ and DS$^o$ into two, but neither of the two partitions bears a relevant meaning. The definition of roaming is the sole purpose of this surface in this investigation.
  
  We can formulate dissociation as a problem of transport of trajectories between dividing surfaces in the following manner. By definition, trajectories that leave the well must cross the outward hemisphere of DS$^i$. Trajectories that leave the interaction region (dissociate), cross the outward annulus of DS$^o$. Dissociation is therefore a question of how trajectories from the outward hemisphere of DS$^i$ reach the outward annulus of DS$^o$.
  
  In the process of being transported from DS$^i$ to DS$^o$, dissociating trajectories originating must cross the outward annulus of DS$^a$. Roaming trajectories cross the torus DS$^a$ at least three times, first the outward annulus and then alternately the inward and the outward annulus. It should be noted that trajectories originating in either of the wells have to cross the outward annulus of DS$^a$ before they can reach its inward annulus. Similarly trajectories going from the inward annulus of DS$^a$ to the outward annulus of DS$^o$ must cross the outward annulus of DS$^a$. Therefore trajectories that originate in the well and cross the inward annulus of DS$^a$ before dissociating are roaming trajectories.
  
  The role of phase space structures in reaction dynamics has been explained on many occasions \cite{Wigginsetal01,Uzeretal02,Waalkens04}. Dividing surfaces as defined above can be shown to be the surfaces of minimal flux \cite{Waalkens04} and therefore sit in the narrowest part of the bottleneck, be it a sphere such as DS$^i$ or a torus such as DS$^a$ and DS$^o$. Structures that convey trajectories across such a bottleneck were identified to be stable and unstable invariant manifolds of the unstable periodic orbit. These invariant manifolds have a cylindrical structure and form a barrier between trajectories that are lead through the bottleneck (inside) and those are not (outside).
  
  For each unstable orbit $\Gamma$, there are two branches of stable invariant manifold, i.e. manifolds consisting of a continuum of trajectories that are asymptotic to the orbit in forward time, and two branches of unstable invariant manifold, asymptotic in backward time. We shall denote the invariant manifolds of $\Gamma$ by $W_{\Gamma}$, we distinguish stable and unstable by a superscript, $W_{\Gamma}^s$ and $W_{\Gamma}^u$, and if we refer to an individual branch, we add a $+$ or $-$ to the superscript to indicate, whether the branch leaves the neighbourhood of $\Gamma$ to the $r>r_\Gamma$ side ($+$) or to the  $r<r_\Gamma$ side ($-$).
    
  The precise cylindrical structure of the manifold reflects the geometry of the corresponding periodic orbit in a similar way as it is reflected by the corresponding dividing surface. Manifolds $W_{\Gamma^i}$ are spherical cylinders, whereas $W_{\Gamma^a}$ and $W_{\Gamma^o}$ are toric cylinders. The different geometries can be understood from a phase space perspective.
  
  For any fixed energy, $\Gamma^i$ is a (topological) circle with its center, due to symmetry of the system, on the ray $\theta=0$, $p_r=p_\theta=0$. This property is passed on to all four branches of $W_{\Gamma^i}$, these too are always centered at/symmetric with respect to $\theta=0$, $p_r=p_\theta=0$, no matter how heavily they are deformed by the flow which is smooth.
  
  In contrast, $W_{\Gamma^a}$ and $W_{\Gamma^o}$ are centered at $r=0$, just like the corresponding orbits $\Gamma^a$ and $\Gamma^o$. This property is independent of energy.
    
  All of this is a consequence of the local energy surface geometry that is not uniform throughout the energy surface and can be observed as qualitatively different forms of dynamics - vibration and rotation.

 \section{Phase space structures responsible for roaming}\label{sec:phase roaming}
  We study invariant manifolds and their intersections on a surface of section, namely on an accurate approximation of the outward annulus of DS$^a$. Both annuli of DS$^a$ are transversal to the flow and are bounded by $\Gamma^a$. Suppose the configuration space projection of $\Gamma^a$ can be parametrised using the function $\bar{r}(\theta)$ as $(\bar{r}(\theta), \theta)$ and define $$\rho(r,\theta)=r-\bar{r}(\theta).$$ Clearly $\Gamma^a$ is then given by $$\rho(r,\theta)=0,$$ and the outward annulus of DS$^a$ is defined by all points $(r,\theta,p_r, p_\theta)$ on the energy surface, such that $$\rho(r,\theta)=0,\quad \dot{\rho}(r,\theta)>0.$$ Note that $\dot{\rho}(r,\theta)>0$ is not equivalent to $\dot{r}>0$.
  
\begin{figure}
   \centering
   \includegraphics[width=.9\textwidth]{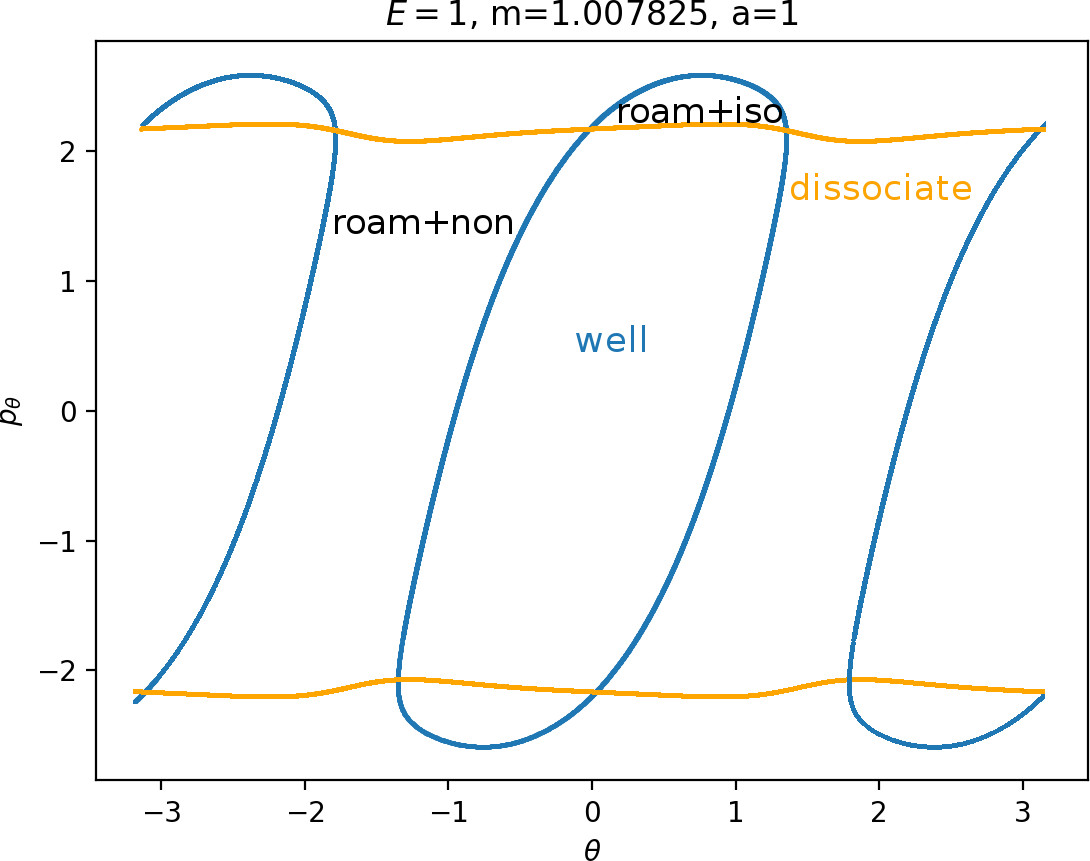}
   \caption{Intersection of $W_{\Gamma^o}^{s-}$ (orange) and $W_{\Gamma^i}^{u+}$ (blue) with the outward annulus of the DS$^a$ for $E=1$, $m=1.007825$ and $a=1$, indicating where trajectories corresponding to different types of dynamics intersect the surface. Area enclosed by $\gamma^{u+}_{i}$ is crossed by trajectories that are leaving the well, while $\gamma^{s-}_{o}$ leads dissociating trajectories out of the interaction region, $\gamma^{u+}_{i}\setminus\gamma^{s-}_{o}$ contains roaming and isomerisation trajectories, and $\gamma^{s-}_{o}\setminus\gamma^{u+}_{i}$ is crossed by roaming and non-reactive trajectories.}
   \label{fig:middlesec}
  \end{figure}

  Most notably, we are interested in the interaction between $W_{\Gamma^i}^{u+}$, the manifold guiding trajectories that leave the potential well into the interaction region, and $W_{\Gamma^o}^{s-}$, the manifold guiding trajectories out of the interaction region into dissociated states. Figure \ref{fig:middlesec} illustrates intersections of these manifolds with the outward annulus of DS$^a$. %for $E=2$ with $a=2,6$ and $m=m_H, 4$.
  
  Note that the section of $W_{\Gamma^i}^{u+}$ with the outward annulus of DS$^a$ produces a topological circle (further denoted $\gamma^{u+}_{i}$) centered at $\theta=0$, $p_\theta=0$. Recall from  Section \ref{sec:phase transport} that all branches of $W_{\Gamma^i}$ are centered at/symmetric with respect to $\theta=0$, $p_r=p_\theta=0$. The shape of $\gamma^{u+}_{i}$ is just a consequence of this fact.
  
  Similarly the section of $W_{\Gamma^o}^{s-}$ with the outward annulus of DS$^a$ (further denoted $\gamma^{s-}_{o}$) reflects the fact that $\Gamma^o$ and all branches of $W_{\Gamma^o}$ are centered at $r=0$. Although $\gamma^{s-}_{o}$ seem like lines, these are two circles concentric with DS$^a$. For simplicity, we omit one of the $\gamma^{s-}_{o}$ circles in all following figure and concentrate on the upper half plane $p_\theta>0$.
  
  All of the figures displaying intersections of invariant manifolds with the outward annulus of DS$^a$ should be understood as follows. Note that for the sake of simplicity we take advantage of symmetry of the system and usually restrict ourselves to $\theta\in[-\pi/2,\pi/2]$, $p_\theta\geq0$ when considering the outward annulus of DS$^a$. Consequently, with the exception of Figure \ref{fig:middlesec}, we only show manifolds corresponding to one of the orbits of $\Gamma^i$ and one of $\Gamma^o$.
  
  Recall that the interior of $\gamma^{u+}_{i}$ contains trajectories that leave the well, while $\gamma^{s-}_{o}$ contains dissociating trajectories which do not return to the surface of section as indicated in Figure \ref{fig:middlesec}. The intersection must therefore contain trajectories that lead to immediate dissociation, such as the radial trajectory $\theta=0$, $p_\theta=0$.
  
  Roaming trajectories do not dissociate immediately, therefore they are contained $\gamma^{u+}_{i}$, but not in $\gamma^{s-}_{o}$. These trajectories have too much energy in the angular degree of freedom, i.e. $|p_\theta|$ is large. The other kind of trajectories contained in $\gamma^{u+}_{i}\setminus\gamma^{s-}_{o}$ are those that re-enter either of the wells, these correspond to isomerisation.
   
  If $\gamma^{u+}_{i}$ and $\gamma^{s-}_{o}$ do not intersect, that is, all of $W_{\Gamma^i}^{u+}$ is contained in $W_{\Gamma^o}^{s-}$, roaming is not present in the system for the given parameter values. For this reason the system does not admit \cite{Krajnak2018} roaming for $E\geq2.5$, $m=m_H$ and $a=1$.
  
  While it is true that $\gamma^{u+}_{i}\setminus\gamma^{s-}_{o}$ contains roaming trajectories, it is not true that the area of the intersection is proportional to the amount of roaming trajectories. Since $\gamma^{u+}_{i}\setminus\gamma^{s-}_{o}$ tends to grow together with $\gamma^{u+}_{i}$ and $\gamma^{u+}_{i}$ can only grow at the cost of $\gamma^{s-}_{o}\setminus\gamma^{u+}_{i}$, the areas $\gamma^{u+}_{i}\setminus\gamma^{s-}_{o}$ and $\gamma^{s-}_{o}\setminus\gamma^{u+}_{i}$ behave like complements. Therefore the amount of roaming is limited by $\gamma^{s-}_{o}\setminus\gamma^{u+}_{i}$, the area where roaming and nonreactive trajectories (see Figure \ref{fig:middlesec}) cross the outer annulus of DS$^a$ for the last time before dissociating.
  
  To see how $m$ and $a$ influence roaming, we will study how $\gamma^{u+}_{i}$ and $\gamma^{s-}_{o}$ change by calculating the areas $\gamma^{u+}_{i}\setminus\gamma^{s-}_{o}$ and $\gamma^{s-}_{o}\setminus\gamma^{u+}_{i}$. The size of $\gamma^{u+}_{i}\setminus\gamma^{s-}_{o}$ is correlated with the maximal value of $p_\theta$ along $\gamma^{u+}_{i}$ and the minimal $p_\theta$ along $\gamma^{s-}_{o}$. Once we know this area, it is easy to determine $\gamma^{s-}_{o}\setminus\gamma^{u+}_{i}$ using the actions of the orbits $\Gamma^i$ and $\Gamma^o$.
  
  The area $\gamma^{u+}_{i}\setminus\gamma^{s-}_{o}$ grows with increasing $m$. Consequently we arrive at one of the main conclusions of this work, namely that roaming diminishes for large and small values of $m$ (Figure \ref{fig:Wiincr}) with an optimum inbetween. This can be seen for $a\leq2$, but for larger values of $a$, the optimum may move towards masses unreasonable for roaming. Disappearance of roaming for large masses is due to a slowdown in the radial degree of freedom and the significant variation of position of $\Gamma^o$ with mass. We explain these two effect separately in Sections \ref{sec:inner} and \ref{sec:outer}. The reason for the disappearance of roaming for low masses is due to a stronger coupling of the degrees of freedom in the kinetic part of the Hamiltonian (\ref{eq:Hamm}) that we deal with in Section \ref{sec:small}.
  
  The parameter $a$ influences the transition between vibration and rotation by controlling the amplitude of the anisotropy in the potential. As noted in Section \ref{sec:potential}, the larger $a$ is, the more the potential wells open up in the angular direction allowing easier access of the wells by trajectories from the interaction region. Instead of promoting roaming, this favours isomerisation, because at the same time the height of the potential barrier between the wells decreases as shown in Figure \ref{fig:thpi2}. We explain the process in Section \ref{sec:alpha}.
  
%   \begin{figure}
%    \centering
%    \includegraphics[width=.49\textwidth]{figure4a.jpg}
%    \includegraphics[width=.49\textwidth]{figure4b.jpg}\\
%    \includegraphics[width=.49\textwidth]{figure4c.jpg}
%    \includegraphics[width=.49\textwidth]{figure4d.jpg}
%    \caption{Intersection of $W_{\Gamma^o_+}^{s-}$ (orange) and $W_{\Gamma^i_+}^{u+}$ (blue) with the outward annulus of the DS$^a$ for $E=2$, masses $m=1.007825,4$ and $a=2,6$.}
%    \label{fig:middlesec}
%   \end{figure}

 \section{Simplification of the system and inner orbits}\label{sec:inner} 
  The mass of the free atom $m$ is only present in terms of the kinetic part of the Hamiltonian (\ref{eq:Hamm}) and therefore cannot be studied from a configuration space perspective. From the Hamiltonian equations of motion
  \begin{eqnarray*}
   \dot{r} &=& \frac{1}{\mu_m} p_r,\\
   \dot{p}_r &=& \frac{p_\theta^2}{\mu_mr^3}-\frac{\partial U}{\partial r},\\
   \dot{\theta} &=& \left(\frac{1}{\mu_m r^2}+\frac{1}{I}\right)p_\theta, \\
   \dot{p}_\theta &=& -\frac{\partial U}{\partial \theta},
  \end{eqnarray*}
  one can see, that $m$ influences the relation between momenta $p_r$, $p_\theta$ and velocities $\dot{r}$, $\dot{\theta}$ as well as the centrifugal contribution in $\dot{p}_r$. All of these are weakened with increasing $m$.
  
  Provided $m$ is sufficiently large (in practice $m>4$), $\frac{1}{\mu r^2}$ is small compared to $\frac{1}{I}$ in the interaction region. For the purposes of the following argument we can neglect the term $\frac{1}{\mu r^2}$ and consider the Hamiltonian
  $$\widetilde{H}_m(r,\theta,p_r, p_\theta) = \frac{1}{2 \mu_m} p_r^2 + \frac{1}{2I}p_\theta^2 + U(r,\theta).$$
  The associated equations are
  \begin{eqnarray*}
   \dot{r} &=& \frac{1}{\mu_m} p_r,\\
   \dot{p}_r &=& -\frac{\partial U}{\partial r},\\
   \dot{\theta} &=& \frac{1}{I} p_\theta, \\
   \dot{p}_\theta &=& -\frac{\partial U}{\partial \theta}.
  \end{eqnarray*}
  The degrees of freedom in the resulting system are only coupled via $U(r,\theta)$ and the mass parameter $m$ only influences $\dot{r}$.
  
  Let $m_0<m_1$, $H_{m_0}=H_{m_1}=E$ and $p_\theta$ be such that it satisfies the fixed energy constraint. Denote by $p_{r_{m_0}}$, $p_{r_{m_1}}$ the momenta given by the Hamiltonians defined as follows:
  \begin{eqnarray*}
   p_{r_{m_0}}=\sqrt{2 \mu_{m_0}\left(E-U(r,\theta)-\frac{p_\theta^2}{2I}\right)},\\
   p_{r_{m_1}}=\sqrt{2 \mu_{m_1}\left(E-U(r,\theta)-\frac{p_\theta^2}{2I}\right)}.
  \end{eqnarray*}
  Then provided $p_\theta$ is not maximal, by
  $$\frac{\mu_{m_0}}{\mu_{m_1}}=\frac{\frac{m_{CH_3}m_0}{m_{CH_3}+m_0}}{\frac{m_{CH_3}m_1}{m_{CH_3}+m_1}}=\frac{m_0m_{CH_3}+m_0m_1}{m_1m_{CH_3}+m_0m_1}<1,$$
  we have that
  \begin{equation*}
   \frac{p_{r_{m_1}}}{\mu_{m_1}}=\frac{1}{\mu_{m_1}}\sqrt{2 \mu_{m_1}\left(E-U-\frac{p_\theta^2}{2I}\right)}
   =\sqrt{\frac{\mu_{m_0}}{\mu_{m_1}}}\frac{1}{\mu_{m_0}}\sqrt{2 \mu_{m_0}\left(E-U-\frac{p_\theta^2}{2I}\right)}
   =\sqrt{\frac{\mu_{m_0}}{\mu_{m_1}}}\frac{p_{r_{m_0}}}{\mu_{m_0}}<\frac{p_{r_{m_0}}}{\mu_{m_0}}.
  \end{equation*}
  For $p_\theta$ maximal we trivially have $p_{r_{m_1}}=p_{r_{m_0}}=0$.
  Trajectories therefore slow down in the radial direction with increasing $m$ for each $r, \theta, p_\theta$ on the energy surface.
  
  Implications for the section of $W_{\Gamma^i}^{u+}$ with the DS$^a$:
  \begin{itemize}
   \item Slowdown in radial direction, whereby $W_{\Gamma^i}^{u+}$ is more likely to hit the potential island.
   \item The interaction of $W_{\Gamma^i}^{u+}$ with the potential island is the only mechanism to transfer energy between the degrees of freedom.
   \item A push from the potential island in the radial direction means decrease in $p_\theta$ along $W_{\Gamma^i}^{u+}$.
%    \item Points with largest $\theta$ correspond to part of manifold passing over the $\widetilde{q}_1^+$ saddle, meaning lowest push in $p_r$ and longest time for rotation.
  \end{itemize}
  
\begin{figure}
   \centering
   \includegraphics[width=.49\textwidth]{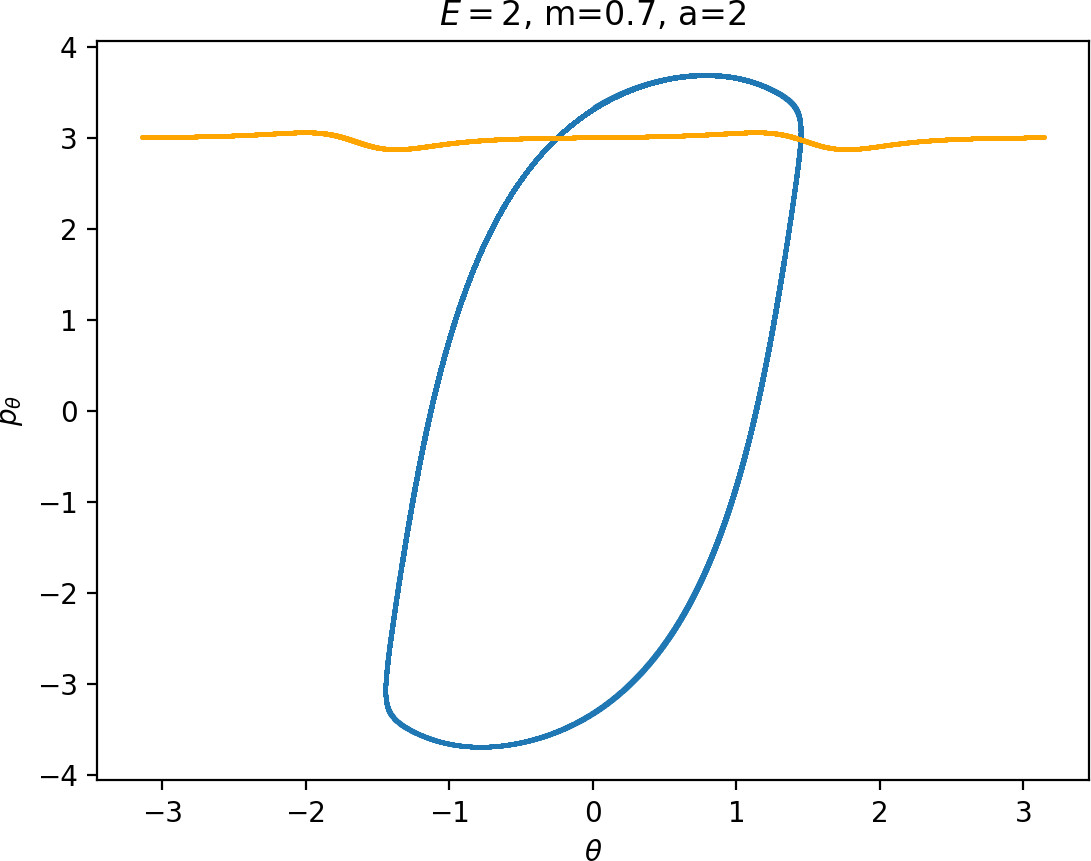}
   \includegraphics[width=.49\textwidth]{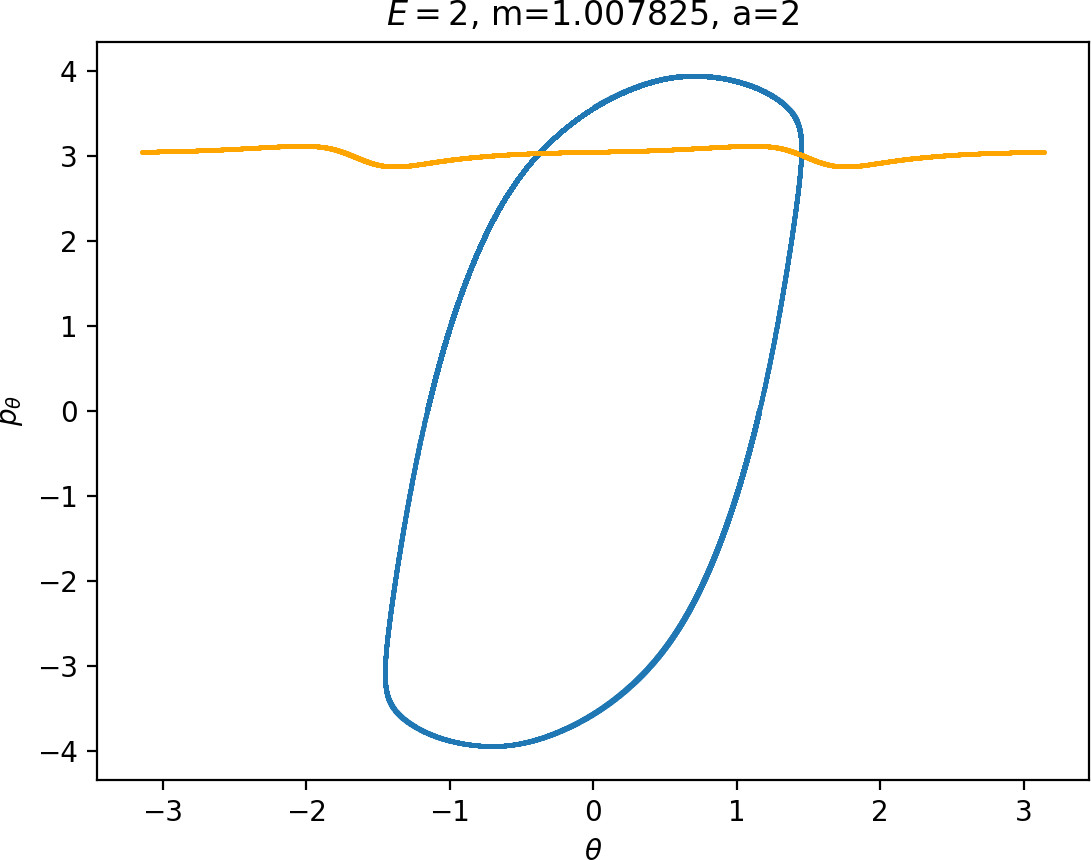}\\
   \includegraphics[width=.49\textwidth]{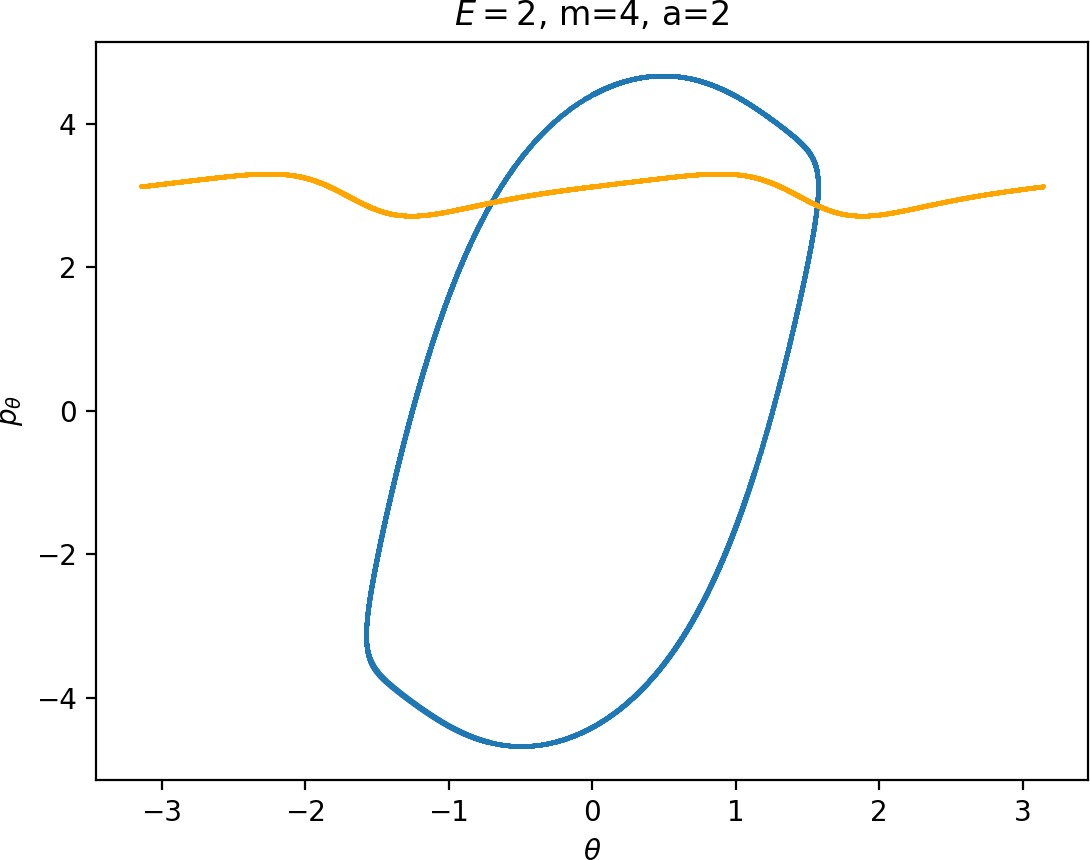}
   \includegraphics[width=.49\textwidth]{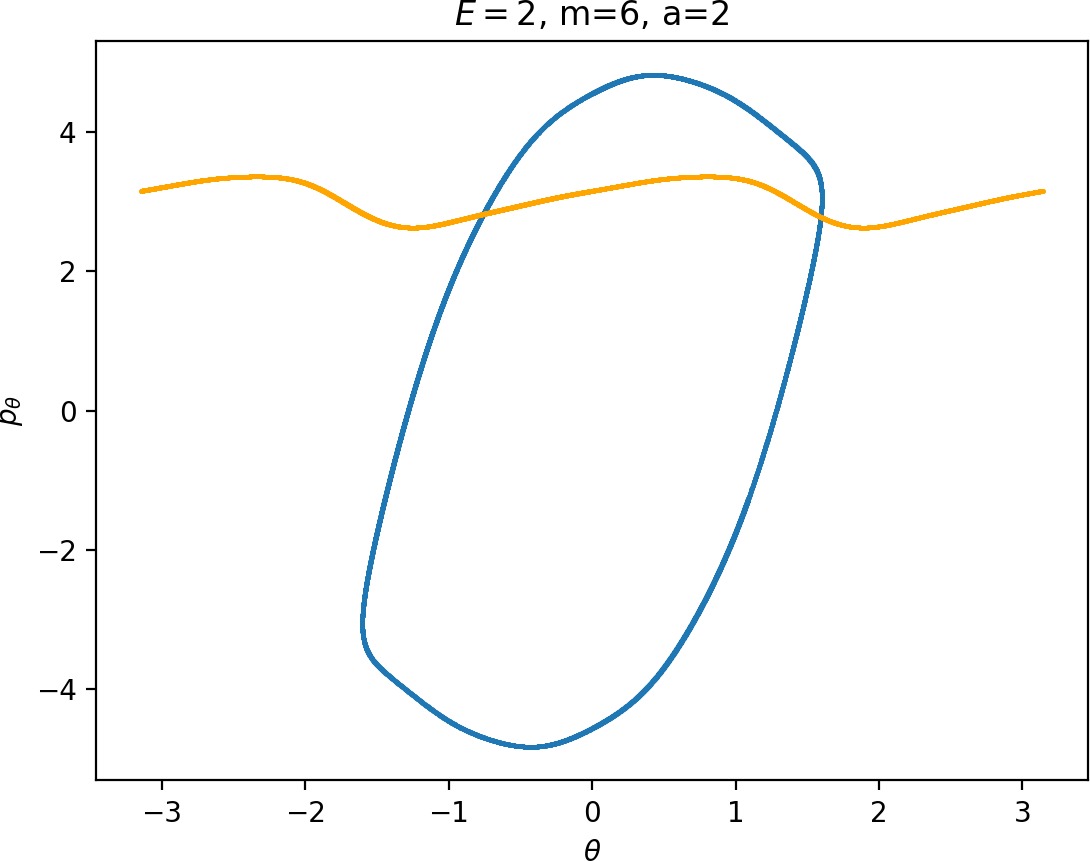}
   \caption{Intersection of $W_{\Gamma^o}^{s-}$ (orange) and $W_{\Gamma^i}^{u+}$ (blue) with the outward annulus of the DS$^a$ for $E=2$, masses $m=0.7,1.007825,4,6$ and $a=2$.}
   \label{fig:Wiincr}
  \end{figure}

\begin{figure}
   \centering
   \includegraphics[width=.49\textwidth]{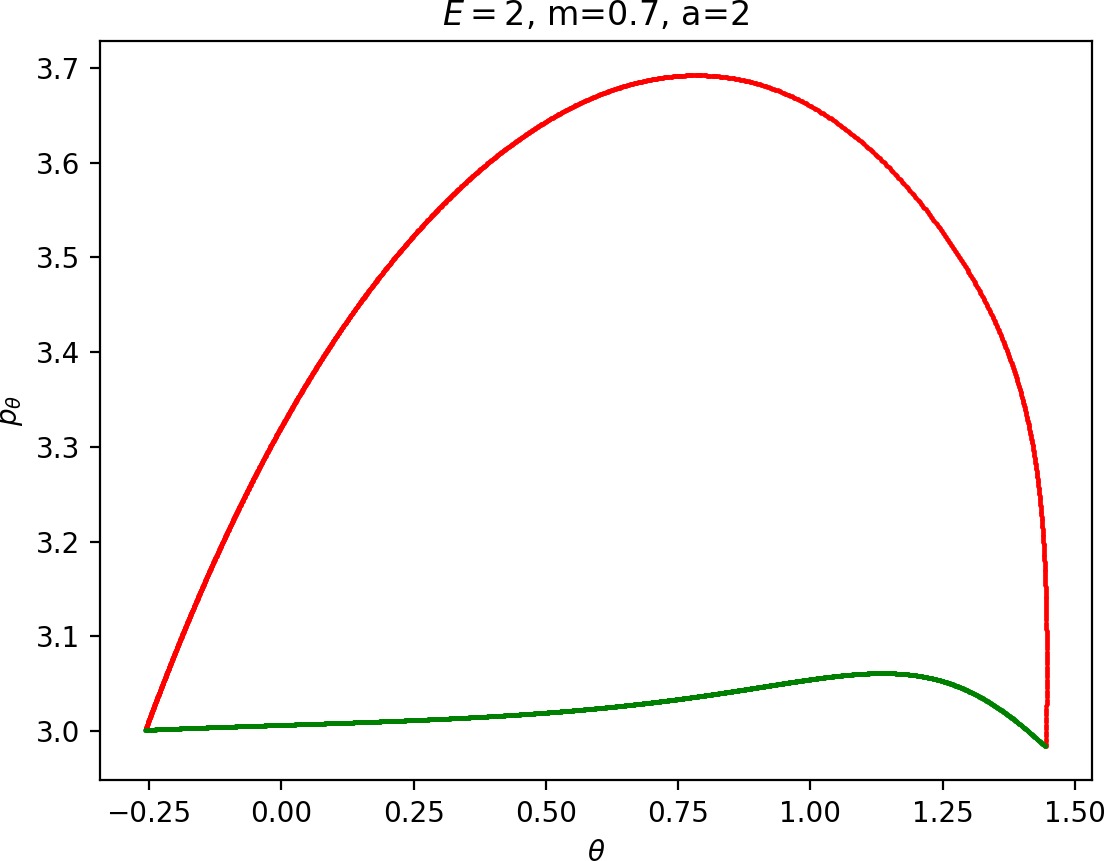}
   \includegraphics[width=.49\textwidth]{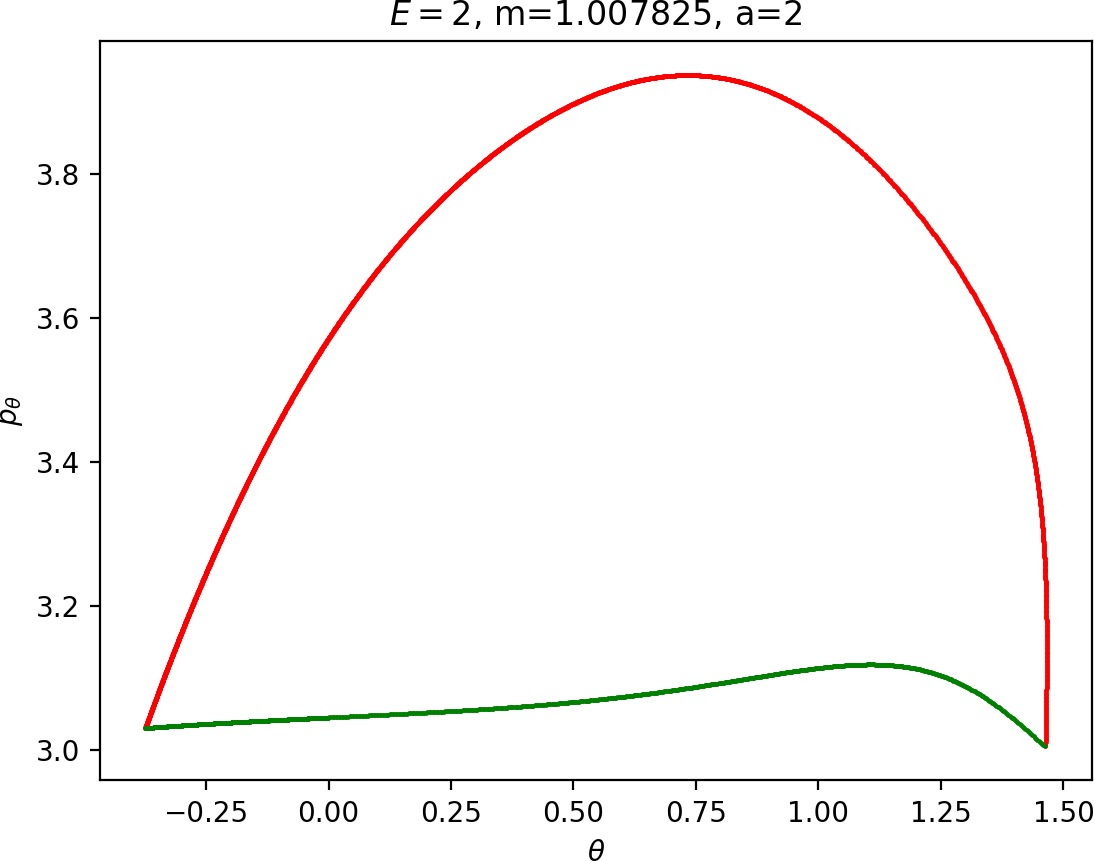}\\
   \includegraphics[width=.49\textwidth]{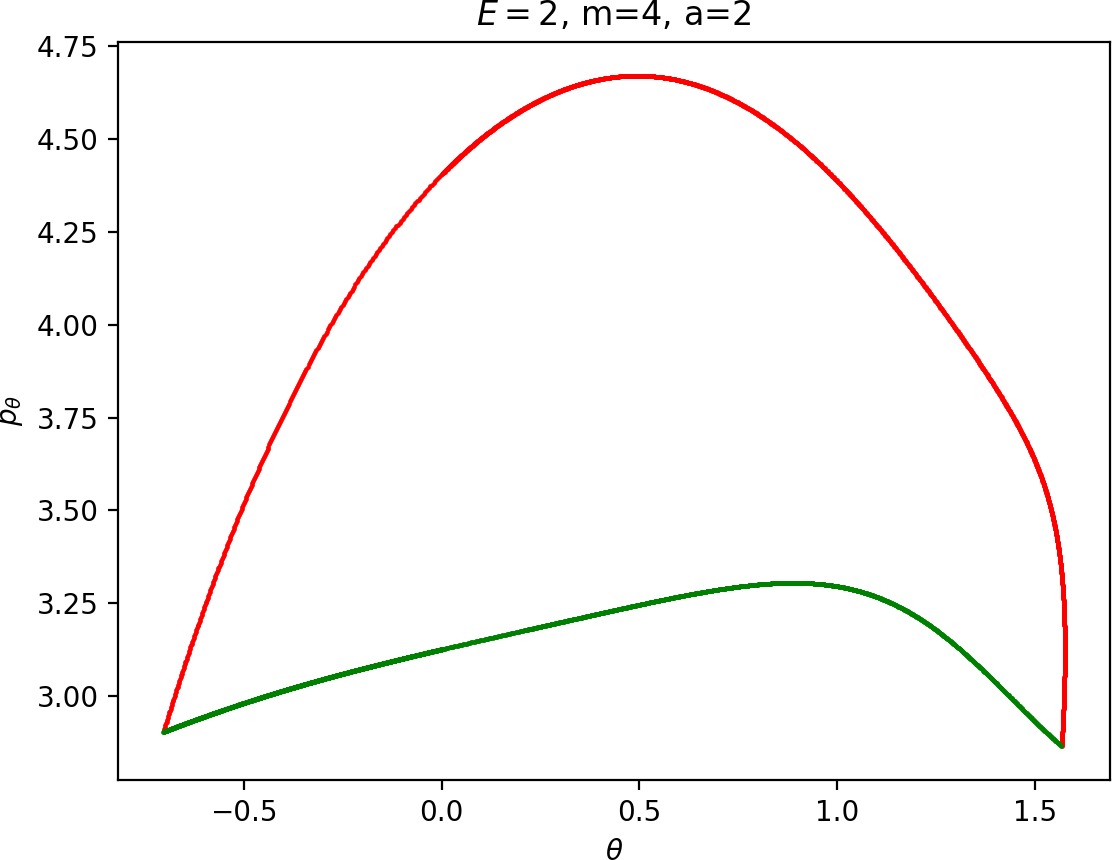}
   \includegraphics[width=.49\textwidth]{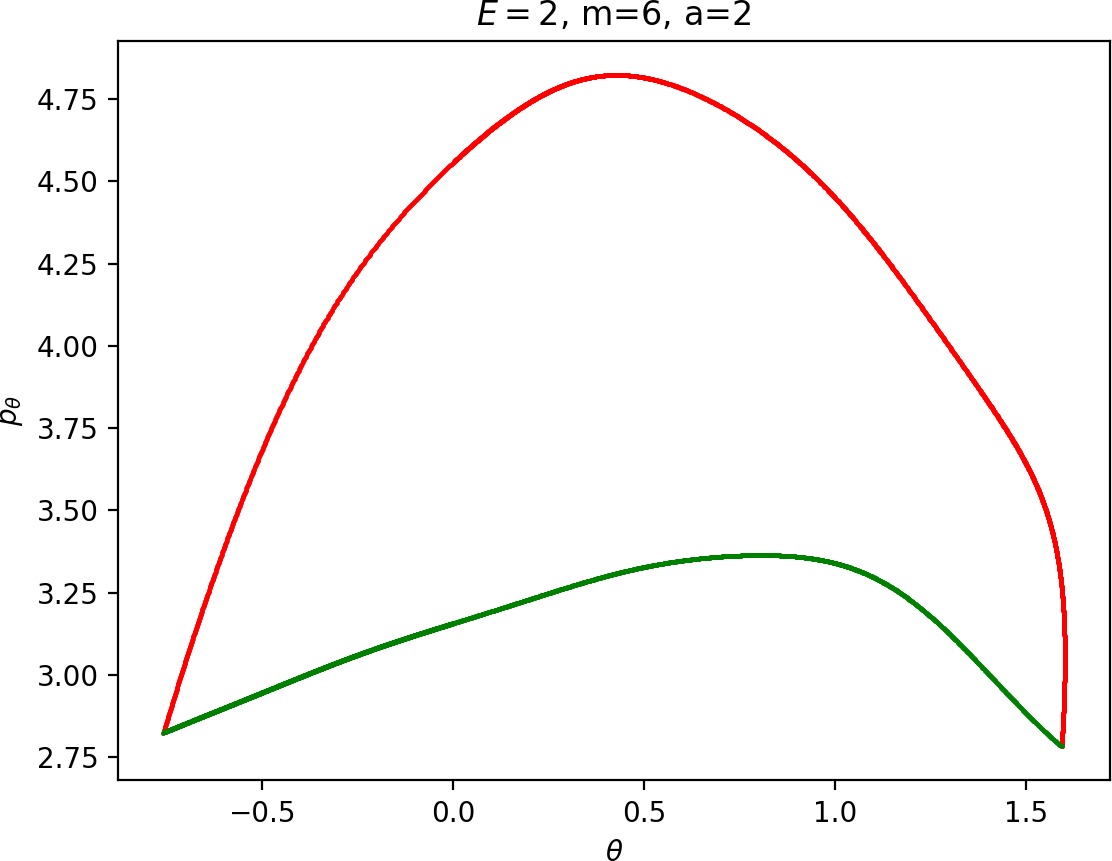}
   \caption{Detail of the intersection of $W_{\Gamma^o}^{s-}$ (green) and $W_{\Gamma^i}^{u+}$ (red) with the outward annulus of the DS$^a$ for $E=2$, masses $0.7,1.007825,4,6$ and $a=2$.}
   \label{fig:Wiincrdetail}
  \end{figure}

  Note that Figures \ref{fig:Wiincr} and \ref{fig:Wiincrdetail} seem to contradict the last point. This is due to the varying radial position of DS$^a$, when $\Gamma^a$ comes closer to $\Gamma^i$ as $m$ increases (see Figure \ref{fig:poconfig}). The decrease can be observed on surfaces with constant radius, however $W_{\Gamma^i}^{u+}$ can become tangent to such surfaces when $m$ is varied.
  
 \section{Reduction of the system and outer orbits}\label{sec:outer}
  If $r$ is sufficiently large, $U$ is rotationally symmetric and
  $$V_{red}(r)=\frac{p_\theta^2}{2 m r^2} + U(r),$$
  where $p_\theta$ becomes a constant of motion. $U(r)$ is monotonic, $U(r)<0$ and $U\in o(r^{-2})$ as $r\rightarrow\infty$, therefore the reduced system admits an equilibrium given by $r=r_{p_\theta}$, $p_r=0$, where $r_{po}$ is the solution of
  \begin{equation*}
  \dot {p}_r = -\frac{\partial H}{\partial r} = \frac{1}{\mu r_{po}^3} p_\theta^2 - U'(r_{po}) = 0, 
  \end{equation*}
  or equivalently
  \begin{equation}\label{eq:pr}
  \frac{p_\theta^2}{\mu r_{po}^2} =r_{po} U'(r_{po}).
  \end{equation}
  In the full system this relative equilibrium is manifested as a periodic orbit with initial conditions $r=r_{po}$, $\theta=const$, $p_r=0$ and $p^\pm_\theta$ such that $H(r_{po}, 0, \theta, p^\pm_\theta)=E$.
  
  To understand the influence of $m$ on the periodic orbit, let us first establish the relationship between $r_{po}$ and $p^\pm_\theta$ and incorporate the dependence on $\mu_m$ subsequently. Combining $H(r_{po}, 0, \theta, p^\pm_\theta)=E$, i.e.
  \begin{equation}\label{eq:Epr0}
   E=\frac{{p^\pm_\theta}^2}{2I}+\frac{{p^\pm_\theta}^2}{2\mu_m r_{po}^2}+U(r_{po}),
  \end{equation}  
  with (\ref{eq:pr}) yields
  \begin{equation*}%\label{eq:tso}
   E=\frac{{p^\pm_\theta}^2}{2I}+\frac{1}{2}r_{po} U'(r_{po})+U(r_{po}).
  \end{equation*}
  If $\frac{1}{2}r U'(r)+U(r)$ is monotonic in $r$ (for $r$ sufficiently large), then we find that $r_{po}$ increases or decreases with ${p^\pm_\theta}^2$. In case of Chesnavich's CH$_4^+$ model, $\frac{1}{2}r U'(r)+U(r)$ is positive and monotonically decreasing. This is due to the leading term of $U$ for large $r$ being $-cr^{-4}$, where $c>0$, while that of $\frac{1}{2}r U'(r)$ must be $2cr^{-4}$. It is sufficient that $U(r)<0$ and $U\in o(r^{-2})$ as $r\rightarrow\infty$. Consequently we see that an increase in ${p^\pm_\theta}^2$ must lead to an increase in $r_{po}$ and vice versa for every fixed energy $E$.
  
  To gain insight on the influence of $m$, we rewrite (\ref{eq:pr}) as
  $${p^\pm_\theta}^2 =\mu_m r_{po}^3 U'(r_{po}),$$
  and plug it into (\ref{eq:Epr0}) to obtain
  \begin{equation*}
   E=\frac{\mu_m r_{po}^3 U'(r_{po})}{2I}+\frac{\mu_m r_{po}^3 U'(r_{po})}{2\mu_m r_{po}^2}+U(r_{po})=\frac{\mu_m r_{po}^3 U'(r_{po})}{2I}+\frac{1}{2}r_{po}U'(r_{po})+U(r_{po}).
  \end{equation*}
  We have already noted that (for $r$ sufficiently large) $\frac{1}{2}r U'(r)+U(r)$ is positive and monotonically decreasing. As a function of $r$ (where $r$ is sufficiently large), $\frac{\mu_m r^3 U'(r)}{2I}$ is also positive and monotonically decreasing, provided $U\in o(r^{-3})$ which is the case for Chesnavich's potential. Therefore to maintain a fixed energy $E$ and increase in $\mu_m$ must be compensated by an increase in $r_{po}$.
  
  We see that $r_{po}$ and ${p^\pm_\theta}^2$ increase with $m$. Note that the reasoning remains true any potential that is rotationally symmetric, monotonic and in the class $o(r^{-3})$ for $r$ sufficiently large and can be extended to systems with non-zero total angular momentum.
  
  Implications for $\gamma^{s-}_{o}$:
  \begin{itemize}
   \item $\Gamma^o_\pm$ moves away with increasing $m$.
   \item $p_\theta$ increases along $\gamma^{s-}_{o}$ with $m$ and it is less influenced by the radial degree of freedom.
  \end{itemize}
   
 \section{Small masses}\label{sec:small}
\begin{figure}
   \centering
   \includegraphics[width=.49\textwidth]{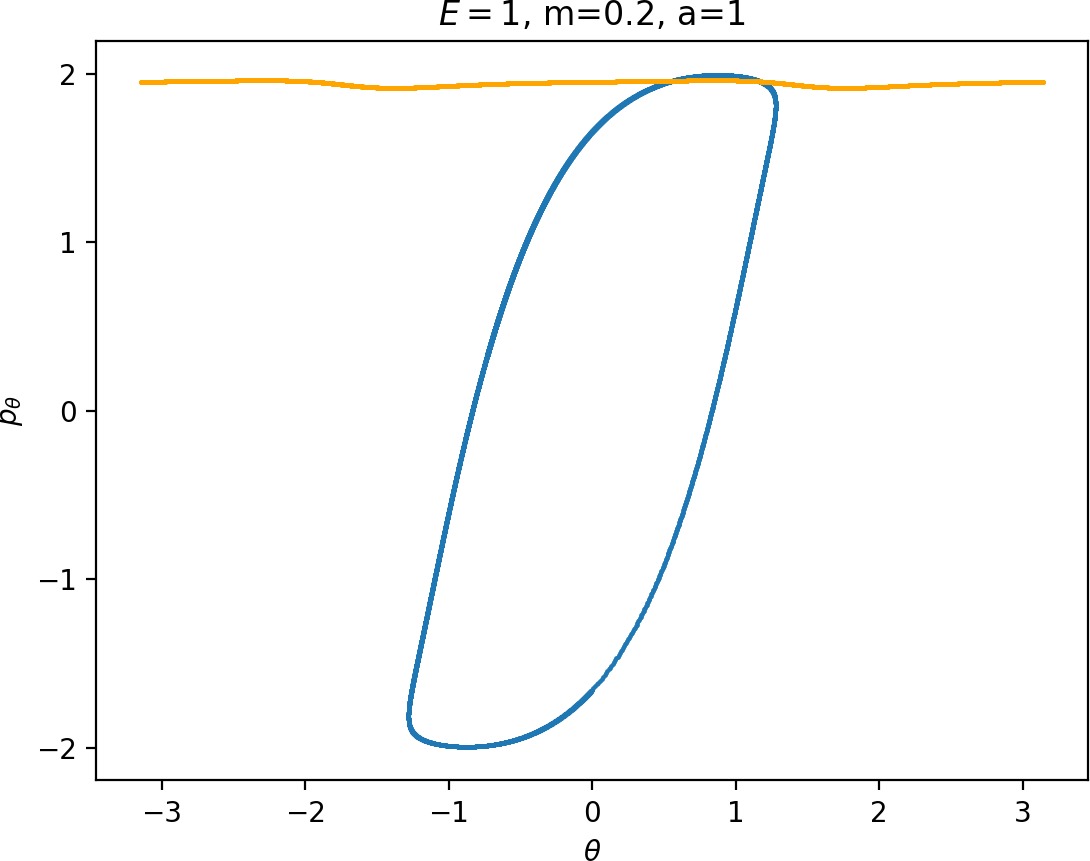}
   \includegraphics[width=.49\textwidth]{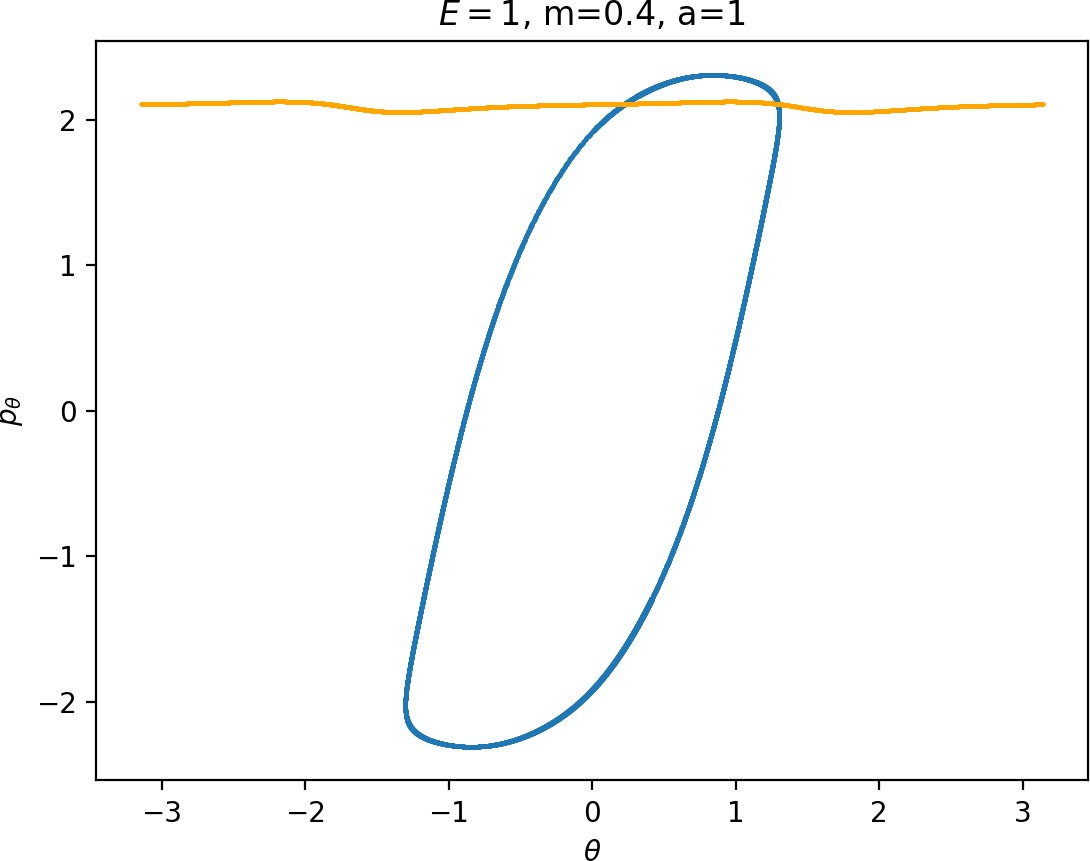}\\
   \includegraphics[width=.49\textwidth]{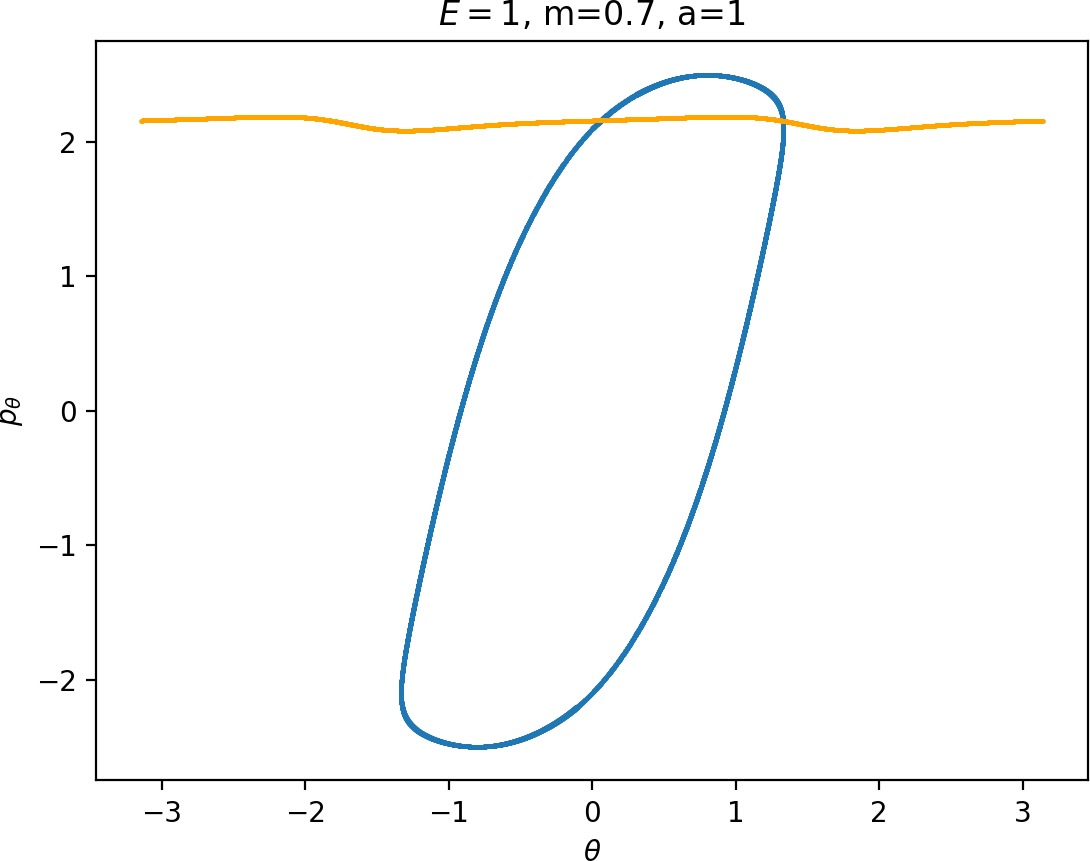}
   \includegraphics[width=.49\textwidth]{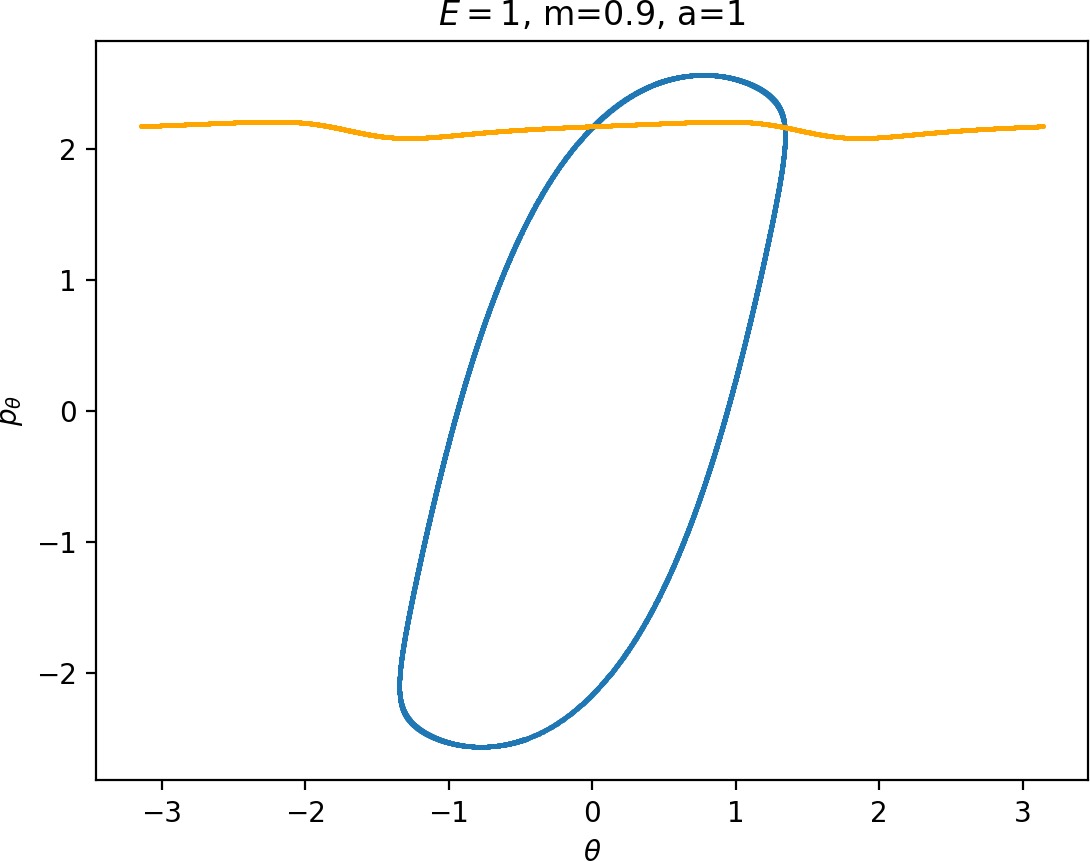}
   \caption{Intersection of $W_{\Gamma^o}^{s-}$ (orange) and $W_{\Gamma^i}^{u+}$ (blue) with the outward annulus of DS$^a$ for $E=1$, $a=1$ and masses $m=0.2,0.4,0.7,0.9$.}
   \label{fig:middlesmallm}
  \end{figure}

  As justified in Section \ref{sec:outer}, $r_{po}$ and ${p^\pm_\theta}^2$ increase with $m$. Therefore for $m<m_H$ we see that $\Gamma^o_\pm$ moves inward and $p_\theta$ along $\gamma^{s-}_{o}$ decreases. This alone does not suffice to make conclusions regarding roaming, because we know that roaming does not exist for atoms with very small masses as suggested by Figure \ref{fig:middlesmallm}.
  
  For $m<m_H$, $p_\theta^2$ has a significant influence on $\dot{p}_r$, especially in the wells where $r$ is small, and also $\dot{r} = \frac{1}{m} p_r$ grows faster. This prevents $W_{\Gamma^i}^{u+}$ from entering areas of high potential near the potential islands and limits the transfer of energy into the angular degree of freedom, since the contribution of the potential in the interaction region is minimal. Although $p_\theta$ along $\Gamma^o$ decreases, %whereby the phase space bottleneck becomes narrower,
  fewer and shorter segments of $W_{\Gamma^i}^{u+}$ can gain sufficient angular momentum to be repelled by the centrifugal barrier back into the interaction region.

 \section{Influence of $a$ on roaming}\label{sec:alpha}
  The parameter $a$ influences the strength of coupling of the degrees of freedom via the potential and thereby it influences \cite{Chesnavich1986} how early (e.g. $a=4$) or late (e.g. $a=1$) the transition from vibration to rotation occurs. An investigation of the differences in dynamics for $a=1$ and $a=4$ can be found in Ref. \onlinecite{Mauguiere2014b}. In physical terms, $a$ controls the anisotropy of the rigid molecule $CH_3^+$.
  
\begin{figure}
   \centering
   \includegraphics[width=.49\textwidth]{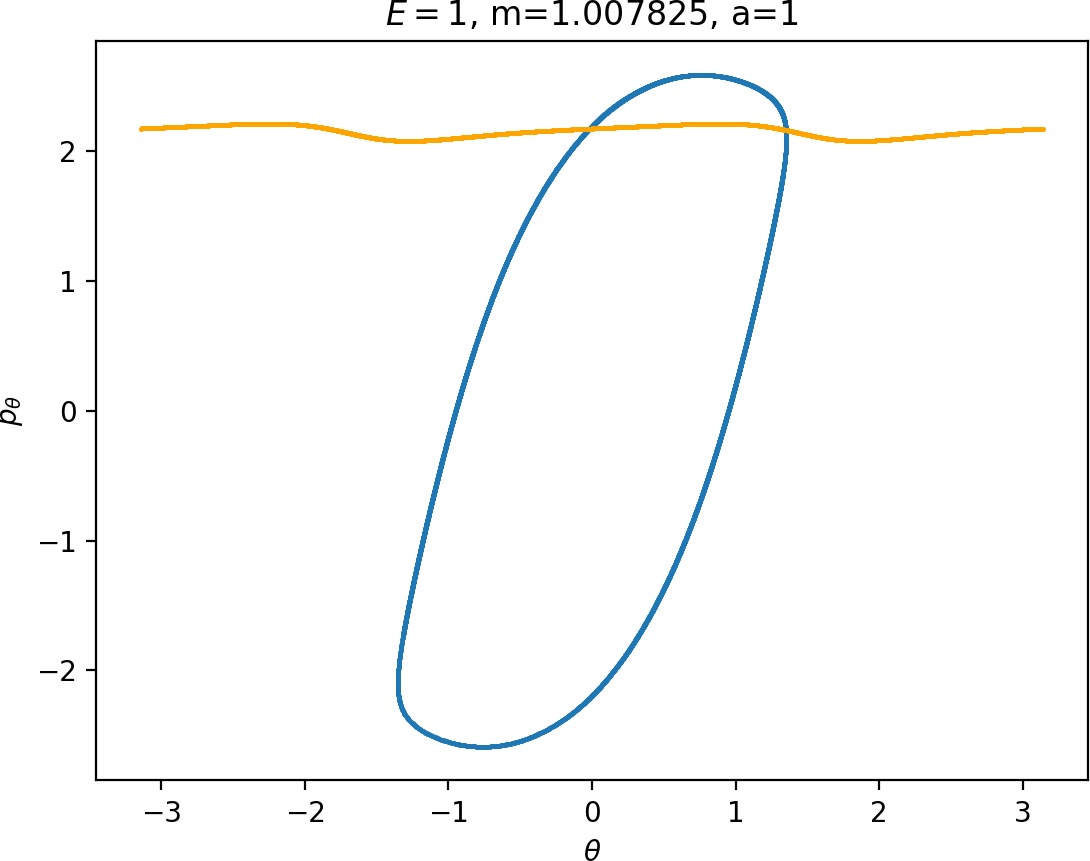}
   \includegraphics[width=.49\textwidth]{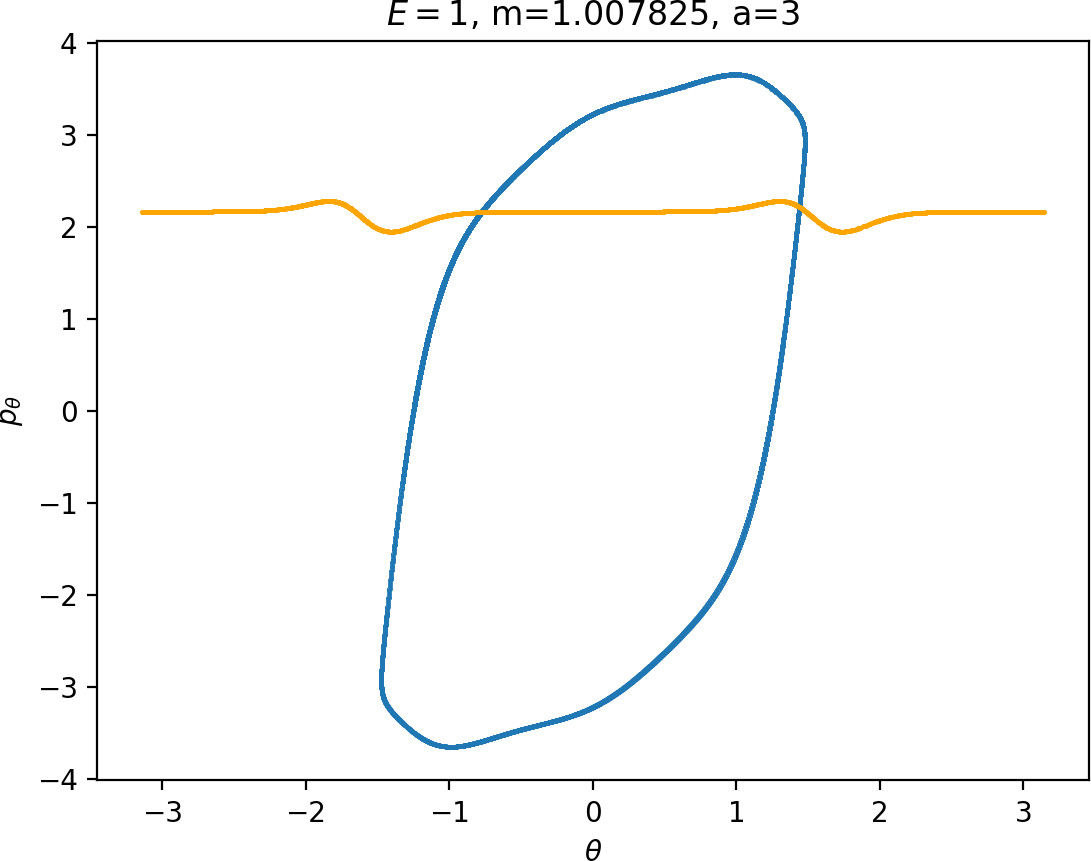}\\
   \includegraphics[width=.49\textwidth]{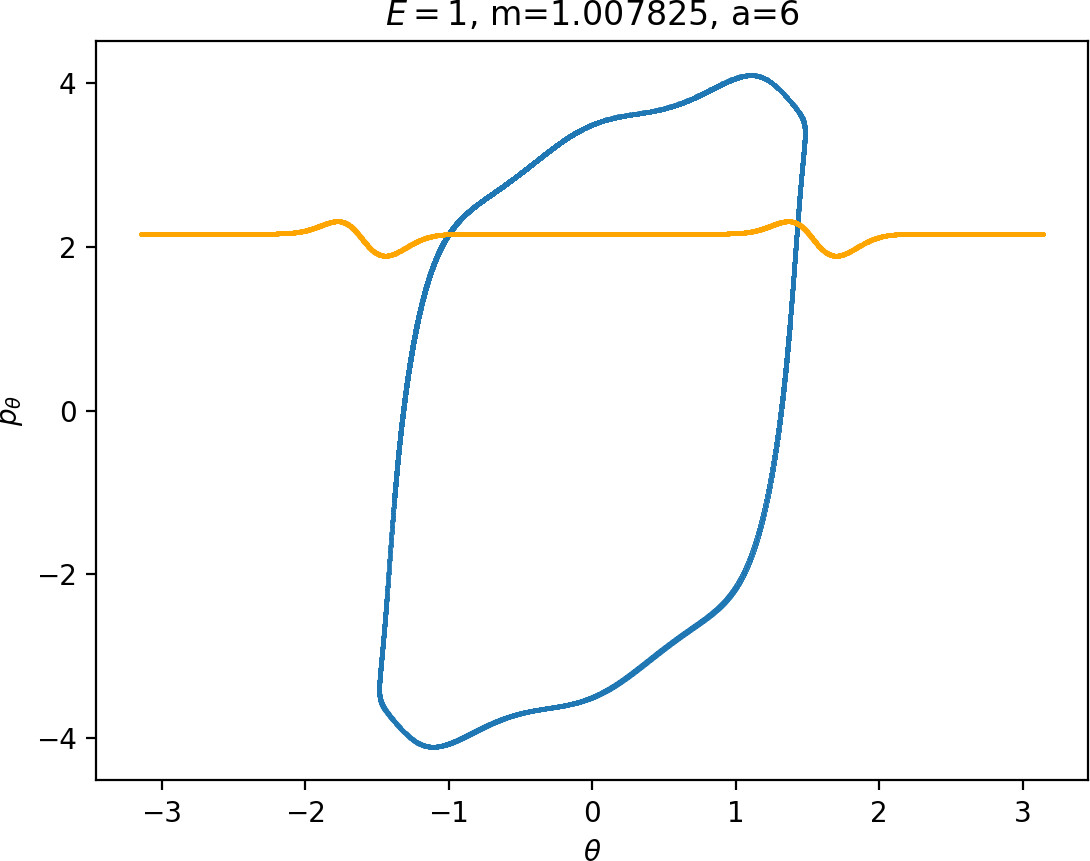}
   \includegraphics[width=.49\textwidth]{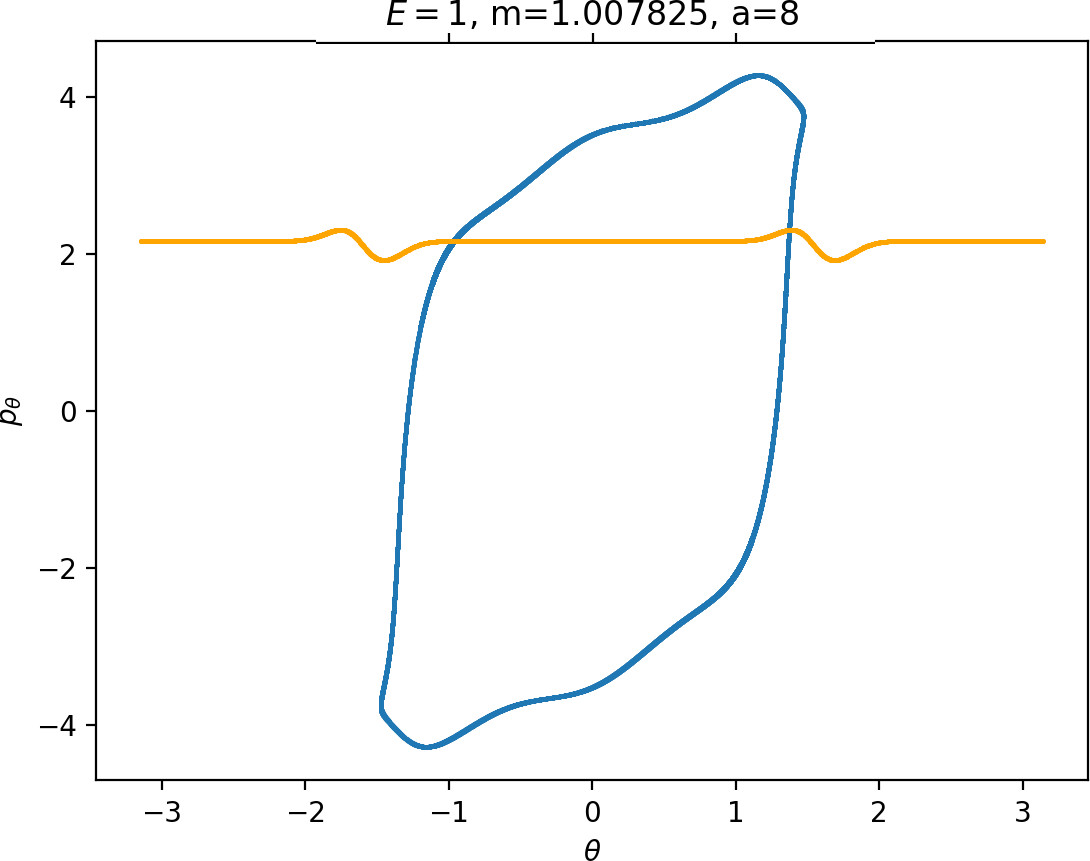}
   \caption{Intersection of $W_{\Gamma^o}^{s-}$ (orange) and $W_{\Gamma^i}^{u+}$ (blue) with the outward annulus of the DS$^a$ for $E=1$, mass $m_H$ and coupling $a=1,3,6,8$.}
   \label{fig:middleE1aa}
  \end{figure}

  It turns out that the expansion of the potential wells and the reduction of potential islands around the index-$2$ saddles with increasing $a$ does have a significant impact on roaming. As can be seen from Figure \ref{fig:middleE1aa}, the area $\gamma^{u+}_{i}\setminus\gamma^{s-}_{o}$ increases with $a$. However, as pointed out in Section \ref{sec:phase roaming} if $\gamma^{u+}_{i}\setminus\gamma^{s-}_{o}$ grows, $\gamma^{s-}_{o}\setminus\gamma^{u+}_{i}$ must shrink and therefore larger values of $a$ result all other classes of dynamics diminishing in favour of isomerisation.
  
 \section{Resulting upper bound on roaming}
\begin{figure}
   \centering
   \includegraphics[height=.3\textheight]{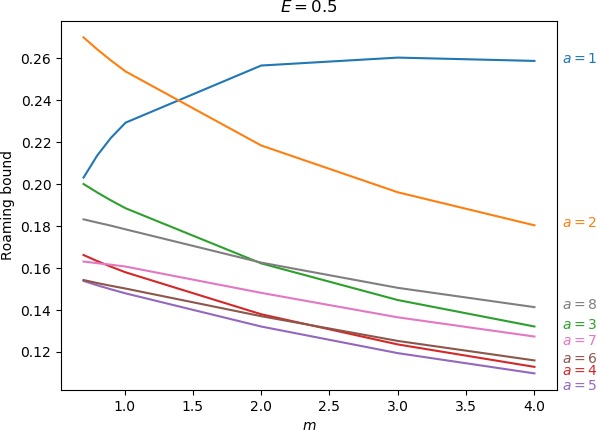}\\
   \includegraphics[height=.3\textheight]{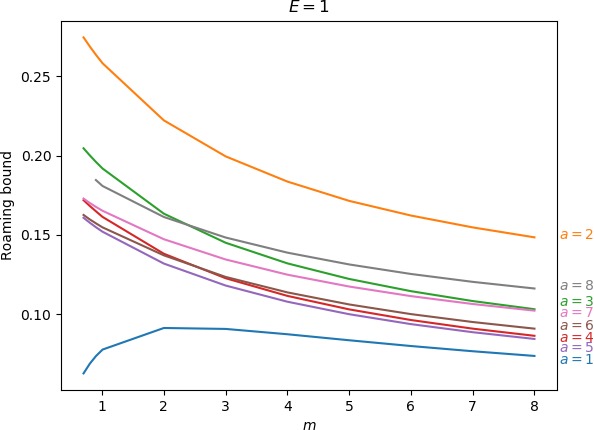}\\
   \includegraphics[height=.3\textheight]{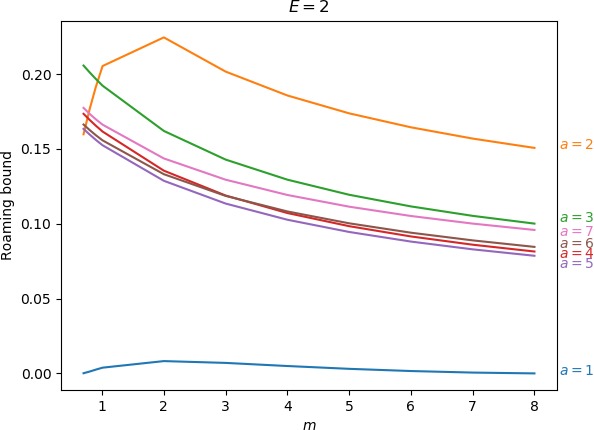}
   \caption{Upper bound on roaming for energies $E=0.5,1,2$.}
   \label{fig:corrroam}
  \end{figure}

  Below we give the upper bound on roaming, which is the ratio of the minimum of the areas $\gamma^{u+}_{i}\setminus\gamma^{s-}_{o}$ and $\gamma^{s-}_{o}\setminus\gamma^{u+}_{i}$ to the energy surface volume for $r=\infty$. As explained above, the area $\gamma^{u+}_{i}\setminus\gamma^{s-}_{o}$ contains all isomerisation and roaming trajectories, while $\gamma^{s-}_{o}\setminus\gamma^{u+}_{i}$ contains all roaming and nonreactive trajectories. The energy surface volume at $r=\infty$ corresponds to all trajectories in the system with the exception of stable islands, which are not relevant in the context of roaming. This measure does not account for isomerisation trajectories, because these are cut off when they re-enter one of the wells. Since all trajectories with the exception of stable islands must eventually leave the wells and the interaction region, we effectively avoid double-counting. Values of the upper bound can be found in Table \ref{tab:1} as well as in Figure \ref{fig:corrroam}.
  
  We do not provide values for $E=0.5$ and $m>4$ due to the large radius of $\Gamma^o$; for $E=1$, $a=8$ and $m<0.9$ due to a bifurcation of $\Gamma^i$ around $m=0.856$ after which DS$^i$ no longer delimits the corresponding potential well in a reasonable manner and similarly for $E=2$ and $a=8$.
  
  As the values in Table \ref{tab:1} show, the prominence of roaming does not evolve in a simple manner. In agreement with conclusions in Sections \ref{sec:inner}, \ref{sec:outer}, \ref{sec:small}, the proportion of roaming decreases or even disappears for very large and very small masses of the free atom and we known it also recedes in favour of isomerisation with increasing $a$. Unlike the area of $\gamma^{u+}_{i}\setminus\gamma^{s-}_{o}$, which grows monotonically both in $m$ and $a$, our upper bound on roaming being the smaller of $\gamma^{u+}_{i}\setminus\gamma^{s-}_{o}$ and $\gamma^{s-}_{o}\setminus\gamma^{u+}_{i}$ does not. Figure \ref{fig:corrroam} suggest that while there is an optimum at a reasonable mass for $a=1$ and occasionally $a=2$, the bound otherwise decreases with $m$ with an optimum at a value of $m$ smaller than considered in this work.
  
  The strange evolution of the bound with respect to $a$, which is shown in Figure \ref{fig:corrroam}, is probably due to the inaccuracy of the bound. Overall it seems that the proportion of roaming mostly decreases, which is due to prevalence of isomerisation as $a$ increases.
  
  An estimate of the amount of roaming in the system can be obtained via an approximation of the proportion of nonreactive trajectories in $\gamma^{s-}_{o}\setminus\gamma^{u+}_{i}$ either via a computationally expensive brute force method or from the intersection of invariant manifolds with a surface $\theta=const$ that does not suffer from transition between spherical and toric local geometries. Such an estimate is outside of the scope of this work.
  
  \begin{table}
   \begin{tabular}{l l|c c c c c c c c c c c c}
     \multicolumn{2}{c}{$E=0.5$} & \multicolumn{11}{c}{$m$}\\
     & & 0.7 & 0.8 & 0.9 & $m_H$ & 2 & 3 & 4 & 5 & 6 & 7 & 8\\
    \hline
    \multirow{8}{*}{\rotatebox{90}{$a$}}
    & 1 & 0.203 & 0.214 & 0.222 & 0.229 & 0.257 & 0.260 & 0.259 & \hphantom{0.000} & \hphantom{0.000} & \hphantom{0.000} & \hphantom{0.000}\\
    & 2 & 0.270 & 0.264 & 0.259 & 0.254 & 0.218 & 0.196 & 0.180\\
    & 3 & 0.200 & 0.196 & 0.192 & 0.189 & 0.162 & 0.145 & 0.132\\
    & 4 & 0.166 & 0.163 & 0.161 & 0.158 & 0.138 & 0.124 & 0.113\\
    & 5 & 0.154 & 0.152 & 0.150 & 0.148 & 0.132 & 0.119 & 0.110\\
    & 6 & 0.154 & 0.153 & 0.152 & 0.150 & 0.137 & 0.125 & 0.116\\
    & 7 & 0.163 & 0.162 & 0.162 & 0.161 & 0.148 & 0.137 & 0.127\\
    & 8 & 0.183 & 0.182 & 0.180 & 0.178 & 0.163 & 0.151 & 0.141
   \end{tabular}

   \begin{tabular}{l l|c c c c c c c c c c c c}
     \multicolumn{2}{c}{$E=1$\hphantom{$.5$}} & \multicolumn{11}{c}{$m$}\\
     & & 0.7 & 0.8 & 0.9 & $m_H$ & 2 & 3 & 4 & 5 & 6 & 7 & 8\\
    \hline
    \multirow{8}{*}{\rotatebox{90}{$a$}}
      & 1 & 0.063 & 0.069 & 0.074 & 0.078 & 0.091 & 0.091 & 0.087 & 0.084 & 0.080 & 0.077 & 0.074\\
      & 2 & 0.274 & 0.269 & 0.264 & 0.258 & 0.222 & 0.200 & 0.184 & 0.171 & 0.162 & 0.155 & 0.148\\
      & 3 & 0.205 & 0.200 & 0.196 & 0.192 & 0.163 & 0.145 & 0.132 & 0.122 & 0.115 & 0.108 & 0.103\\
      & 4 & 0.172 & 0.168 & 0.165 & 0.161 & 0.138 & 0.123 & 0.112 & 0.103 & 0.096 & 0.091 & 0.086\\
      & 5 & 0.161 & 0.158 & 0.155 & 0.152 & 0.132 & 0.118 & 0.108 & 0.100 & 0.094 & 0.089 & 0.084\\
      & 6 & 0.163 & 0.160 & 0.157 & 0.155 & 0.137 & 0.124 & 0.114 & 0.106 & 0.100 & 0.095 & 0.091\\
      & 7 & 0.173 & 0.170 & 0.168 & 0.165 & 0.147 & 0.135 & 0.125 & 0.117 & 0.111 & 0.106 & 0.102\\
      & 8 &  &  & 0.185 & 0.181 & 0.161 & 0.148 & 0.139 & 0.131 & 0.125 & 0.120 & 0.116
   \end{tabular}

   \begin{tabular}{l l|c c c c c c c c c c c c}
     \multicolumn{2}{c}{$E=2$\hphantom{$.5$}} & \multicolumn{11}{c}{$m$}\\
     & & 0.7 & 0.8 & 0.9 & $m_H$ & 2 & 3 & 4 & 5 & 6 & 7 & 8\\
    \hline
    \multirow{7}{*}{\rotatebox{90}{$a$}}
      & 1 & 0.000 & 0.001 & 0.003 & 0.004 & 0.008 & 0.007 & 0.005 & 0.003 & 0.002 & 0.001 & 0.000\\
      & 2 & 0.160 & 0.177 & 0.192 & 0.205 & 0.224 & 0.202 & 0.186 & 0.174 & 0.164 & 0.157 & 0.151\\
      & 3 & 0.206 & 0.201 & 0.197 & 0.192 & 0.162 & 0.143 & 0.129 & 0.119 & 0.112 & 0.105 & 0.100\\
      & 4 & 0.173 & 0.169 & 0.165 & 0.162 & 0.135 & 0.119 & 0.107 & 0.098 & 0.092 & 0.086 & 0.081\\
      & 5 & 0.163 & 0.160 & 0.156 & 0.152 & 0.129 & 0.113 & 0.103 & 0.095 & 0.088 & 0.083 & 0.079\\
      & 6 & 0.166 & 0.163 & 0.159 & 0.156 & 0.133 & 0.119 & 0.108 & 0.100 & 0.094 & 0.089 & 0.085\\
      & 7 & 0.177 & 0.173 & 0.170 & 0.166 & 0.144 & 0.129 & 0.119 & 0.111 & 0.105 & 0.100 & 0.096\\
   \end{tabular}
   \caption{Ratio of the minimum of the areas $\gamma^{u+}_{i}\setminus\gamma^{s-}_{o}$ and $\gamma^{s-}_{o}\setminus\gamma^{u+}_{i}$ to the measure of all trajectories in the system (for details see text) for $E=0.5,1,2$ and various values of $m$ and $a$.}
   \label{tab:1}
  \end{table}
  
  \section{Summary and conclusion}
   We have shown how invariant manifolds that are responsible for roaming evolve with respect to the mass of the free atom $m$ and coupling parameter $a$ that controls the geometry of the potential energy surface. Moreover we provided arguments that follow from the role of parameters in the equations of motion to justify the evolution. These arguments make it also possible to predict behaviour outside of the studied parameter intervals relevant for roaming.
   
   From a quantitative perspective we established an upper bound on the prominence of roaming in Chesnavich's CH$_4^+$ model. The bound only highlights the intricacy of roaming as a type of dynamics on the verge between isomerisation and nonreactivity. It relies on a generous access to the potential wells to allow reactions as well as a prominent area of high potential that aids sufficient transfer of energy from angular motion to radial motion to prevent isomerisation.
   
   We conclude that it is not possible to choose realistic values of $m$ and $a$, such that roaming becomes the dominant form of dissociation such as found in acetaldehyde \cite{Heazlewood08}. Therefore the dominance of roaming must be due to other properties of the system than only mass of the free atom and strength of potential coupling of the degrees of freedom. Our investigation shows that roaming is prone to being dominated by other types of dynamics. Therefore the path to dominance of roaming lies in, as our investigation shows, potential wells that are easily accessible by being 'open' in the angular direction to discourage nonreactivity, yet are separated by a sufficiently high isomerisation barrier.
   
   It is interesting to speculate about the relationship between the roaming mechanism and reactions of floppy molecules, such as HCN and KCN. Both types of reaction are focussed on large-amplitude motions. As we have described, in Chesnavich's model it is the relatively long-range part of the potential that is of interest; that is where dissociation, roaming, and the passage from one inner DS to the other occurs.

  With HCN, for example, the focus is on isomerisation dynamics. In the settting of Chesnavich's model, this means the concern is about the properties of the `inner' region of the potential energy surface. In the study of roaming in Chesnavich's model the inner region of the potential energy surface (interior to the inner DS) is basically left unspecified.

  Nevertheless, as noted in discussions of Ezra and Houston, it is probably true that roaming is just a special case of large-amplitude motions. Consequently, it may well be possible to identify large amplitude eigenstates of floppy molecules, such as HCN or KCN, that correspond to roaming motion, and these may be associated with periodic orbits in the interaction region. This is an interesting topic for further investigation.
  
%-----------------------------------------------------------------
%-----------------------------------------------------------------
\acknowledgments
 We acknowledge the support of EPSRC Grant No.~EP/P021123/1 and ONR Grant No.~N00014-01-1-0769.
 
 Discussions of Greg Ezra and Paul Houston have provided interesting insights into the connection between the roaming mechanism and the reaction of floppy molecules which should provide inspiration for further study. We are grateful that they have shared their insights with us.

% \bibliographystyle{unsrt}
% \bibliography{refs}

\end{document}